\documentclass[12pt]{article}
\usepackage{amssymb, amsmath, amsthm, amscd}
\usepackage[utf8]{inputenc}
\usepackage[all,cmtip]{xy}
\usepackage{enumitem}
\usepackage{setspace}
\usepackage{mathdots}

\usepackage[mathscr]{euscript}
\usepackage[colorlinks=true]{hyperref}
\hypersetup{colorlinks   = true}
\hypersetup{linkcolor=blue}
\usepackage{tikz}

\begin{document}

\mathchardef\mhyphen="2D
\newtheorem{The}{Theorem}[section]
\newtheorem{Lem}[The]{Lemma}
\newtheorem{Prop}[The]{Proposition}
\newtheorem{Cor}[The]{Corollary}
\newtheorem{Rem}[The]{Remark}
\newtheorem{Obs}[The]{Observation}
\newtheorem{SConj}[The]{Standard Conjecture}
\newtheorem{Titre}[The]{\!\!\!\! }
\newtheorem{Conj}[The]{Conjecture}
\newtheorem{Question}[The]{Question}
\newtheorem{Prob}[The]{Problem}
\newtheorem{Def}[The]{Definition}
\newtheorem{Not}[The]{Notation}
\newtheorem{Claim}[The]{Claim}
\newtheorem{Conc}[The]{Conclusion}
\newtheorem{Ex}[The]{Example}
\newtheorem{Fact}[The]{Fact}
\newtheorem{Formula}[The]{Formula}
\newtheorem{Formulae}[The]{Formulae}
\newtheorem{The-Def}[The]{Theorem and Definition}
\newtheorem{Prop-Def}[The]{Proposition and Definition}
\newtheorem{Lem-Def}[The]{Lemma and Definition}
\newtheorem{Cor-Def}[The]{Corollary and Definition}
\newtheorem{Conc-Def}[The]{Conclusion and Definition}
\newtheorem{Terminology}[The]{Note on terminology}
\newcommand{\C}{\mathbb{C}}
\newcommand{\R}{\mathbb{R}}
\newcommand{\N}{\mathbb{N}}
\newcommand{\Z}{\mathbb{Z}}
\newcommand{\Q}{\mathbb{Q}}
\newcommand{\Proj}{\mathbb{P}}
\newcommand{\Rc}{\mathcal{R}}
\newcommand{\Oc}{\mathcal{O}}
\newcommand{\Vc}{\mathcal{V}}
\newcommand{\Id}{\operatorname{Id}}
\newcommand{\pr}{\operatorname{pr}}
\newcommand{\rk}{\operatorname{rk}}
\newcommand{\del}{\partial}
\newcommand{\delbar}{\bar{\partial}}
\newcommand{\Cdot}{{\raisebox{-0.7ex}[0pt][0pt]{\scalebox{2.0}{$\cdot$}}}}
\newcommand\nilm{\Gamma\backslash G}
\newcommand\frg{{\mathfrak g}}
\newcommand{\fg}{\mathfrak g}
\newcommand{\Oh}{\mathcal{O}}
\newcommand{\Kur}{\operatorname{Kur}}
\newcommand\gc{\frg_\mathbb{C}}
\newcommand\hisashi[1]{{\textcolor{red}{#1}}}
\newcommand\dan[1]{{\textcolor{blue}{#1}}}
\newcommand\luis[1]{{\textcolor{orange}{#1}}}

\begin{center}

{\Large\bf Holomorphic $p$-Contact and $s$-Symplectic Line Bundles}

\end{center}

\begin{center}

{\large Kyle Broder and Dan Popovici}

\end{center}

\vspace{1ex}

\noindent{\small{\bf Abstract.} We generalise the notions of scalar-valued holomorphic $p$-contact and $s$-symplectic structures introduced recently on compact complex manifolds by the second-named author jointly with H. Kasuya and L. Ugarte to their analogues with values in a holomorphic line bundle. We then study the resulting holomorphic $p$-contact and $s$-symplectic manifolds which, unlike their scalar counterparts that are never K\"ahler, can even be projective. In particular, we investigate the (lack of) positivity properties of the canonical bundle of these manifolds when it is given a possibly singular Hermitian fibre metric. One of the tools used is a very recent regularisation result for $m$-psh functions obtained jointly by S. Dinew and the second-named author.}

\vspace{1ex}

\section{Introduction}\label{section:introduction} The following generalisations of the classical notions of holomorphic contact and holomorphic symplectic manifolds were introduced in [KPU25a] and further studied in [KPU25b] in the context of scalar-valued (i.e. $\C$-valued) differential forms.

\begin{Def}\label{Def:hol_p-contact_s-symplectic_scalar} Let $X$ be a compact complex manifold with $\mbox{dim}_\C X =n$.

\vspace{1ex}

\noindent $(1)$ ([KPU25, Definition 1.1.]) When $n=2p+1$, $X$ is said to be a {\bf holomorphic $p$-contact manifold} if there exists a holomorphic $p$-contact structure on it, to wit a smooth $(p,\,0)$-form $\Gamma\in C^\infty_{p,\,0}(X,\,\C)$ such that \begin{eqnarray*}(a)\hspace{2ex} \bar\partial\Gamma=0 \hspace{6ex} \mbox{and} \hspace{6ex} (b)\hspace{2ex} \Gamma\wedge\partial\Gamma\neq 0 \hspace{2ex} \mbox{at every point of} \hspace{1ex} X.\end{eqnarray*}

\vspace{1ex}

\noindent $(2)$ ([KPU25, Definition 3.4.]) When $n=2s$, $X$ is said to be a {\bf holomorphic $s$-symplectic manifold} if there exists a holomorphic $s$-symplectic structure on it, to wit a smooth $(s,\,0)$-form $\Omega\in C^\infty_{s,\,0}(X,\,\C)$ such that \begin{eqnarray*}(i)\,\,\bar\partial\Omega = 0; \hspace{3ex} \mbox{and}  \hspace{3ex} (ii)\,\, \Omega\wedge\Omega \neq 0 \hspace{2ex} \mbox{at every point of}\hspace{1ex} X.\end{eqnarray*} 

\end{Def}

It is easily seen that $p$ must be odd, forcing $n\equiv 3$ mod $4$, when $X$ is holomorphic $p$-contact, while $s$ must be even, forcing $n\equiv 0$ mod $4$, when $X$ is holomorphic $s$-symplectic.

\vspace{2ex}

In the present paper, we further generalise these notions to the setting of differential forms with values in a holomorphic line bundle $L$. Specifically, we introduce and study the following

\begin{Def}\label{Def:hol_p-contact_s-symplectic_bundle} Let $X$ be a compact complex manifold with $\mbox{dim}_\C X =n$. 

\vspace{1ex}

\noindent $(1)$ When $n=2p+1=4l+3$, $X$ is said to be a {\bf (locally) holomorphic $p$-contact manifold} if there exist a holomorphic line bundle $L$ and an $L$-valued holomorphic $p$-contact structure on $X$, to wit a form $\Gamma\in C^\infty_{p,\,0}(X,\,L)$ satisfying the following two conditions: \begin{eqnarray}\label{eqn:L_hol-p-contact_Gamma_introd}(a)\,\,\, \bar\partial\Gamma = 0  \hspace{5ex}\mbox{and}\hspace{5ex} (b)\,\,\,\Gamma\wedge\partial\Gamma \neq 0 \hspace{2ex}\mbox{at every point of}\hspace{2ex} X.\end{eqnarray}  

\vspace{1ex}

\noindent $(2)$ When $n=2s$, $X$ is said to be a {\bf (locally) holomorphic $s$-symplectic manifold} if there exist a holomorphic line bundle $L$ and an $L$-valued holomorphic $s$-symplectic structure on $X$, to wit a form $\Omega\in C^\infty_{s,\,0}(X,\,L)$ satisfying the following two conditions: \begin{eqnarray}\label{eqn:L_hol-s-symplectic_Omega_introd}(a)\,\,\, \bar\partial\Omega = 0 \hspace{5ex}\mbox{and}\hspace{5ex} (b)\,\,\,\Omega\wedge\Omega \neq 0 \hspace{2ex}\mbox{at every point of}\hspace{2ex} X.\end{eqnarray}

\end{Def}

Note that $p$ is supposed to be odd in the $p$-contact case, forcing $n\equiv 3$ mod $4$, while $s$ must be even as a consequence of the assumptions (\ref{eqn:L_hol-s-symplectic_Omega_introd}), forcing $n\equiv 0$ mod $4$, in the $s$-symplectic case.

On the other hand, $\partial\Gamma$ is not a globally defined $L$-valued form on $X$ since $\partial$ is not a $(1,\,0)$-connection on $L$, but 
$\Gamma\wedge\partial\Gamma$ is a well-defined $L^2$-valued holomorphic $(n,\,0)$-form on $X$ when $n=2p+1$ with $p$ odd. This can be seen at least in two ways: by viewing $\Gamma$ as a collection $(\Gamma_\alpha)_\alpha$ of $\C$-valued $(p,\,0)$-forms on the open subsets $U_\alpha$ of a family of local trivialisations of $L$ and noting that the locally defined $\C$-valued $(n,\,0)$-forms $\Gamma_\alpha\wedge\partial\Gamma_\alpha$ glue together into a global $L^2$-valued $(n,\,0)$-form on $X$ (cf. (i) of Observation \ref{Obs:Gamma-wedge-sel-Gamma_global}); or by observing that, if $D'_h$ is the $(1,\,0)$-part of the Chern connection $D_h$ of $L$ equipped with any $C^\infty$ Hermitian fibre metric $h$, then $\Gamma\wedge D'_h\Gamma$ is independent of $h$ and equals the form $\Gamma\wedge\partial\Gamma$ defined by the first approach (cf. (ii) of Observation \ref{Obs:Gamma-wedge-sel-Gamma_global}).

\vspace{1ex}

The classical notion of {\bf $L$-valued holomorphic contact structure} (cf. e.g. [LeB95]) on a complex manifold $X$ of odd dimension $n=2p+1$ is defined to be a holomorphic $(1,\,0)$-form $\eta\in H^{1,\,0}_{\bar\partial}(X,\,L)$ such that the $L^{p+1}$-valued holomorphic $(n,\,0)$-form $\eta\wedge(\partial\eta)^p$ is nowhere vanishing. We observe in Proposition \ref{Prop:1-contact_p-contact} that every compact complex manifold $X$ with $\mbox{dim}_\C X = n\equiv 3$ mod $4$ that carries a bundle-valued holomorphic contact structure also carries a bundle-valued holomorphic $p$-contact structure. Thus, our (locally) holomorphic $p$-contact manifolds generalise the classical contact manifolds.

In a similar vein, our locally holomorphic $s$-symplectic manifolds generalise the classical {\bf holomorphic symplectic manifolds} in two ways. Firstly, in the scalar case, it was observed in [KPU25a, $\S3.2$] that every holomorphic symplectic form (i.e. every holomorphic $(2,\,0)$-form $\omega\in H^{2,\,0}_{\bar\partial}(X,\,\C)$ on a $(2s)$-dimensional compact complex manifold $X$ such that $d\omega = 0$ and $\omega^s\neq 0$ everywhere on $X$) induces a $\C$-valued holomorphic $s$-symplectic form defined as $\Omega:=\omega^s$. Secondly, we allow our holomorphic $s$-symplectic forms $\Omega$ to assume values in a possibly non-trivial holomorphic line bundle $L$. 

Meanwhile, the existence of an $L$-valued holomorphic $p$-contact or $s$-symplectic structure implies that the holomorphic line bundle $K_X\otimes L^2$ is trivial, or equivalently that $L$ is a square root of the anti-canonical bundle $-K_X$. This, in turn, is known to be equivalent to the existence of a {\it spin structure} on $X$ (cf. [Ati71], cited here as Proposition \ref{Prop:spin_existence_standard}). Thus, we have

\begin{Obs}\label{Obs:p-contact_s-symplectic_spin_introd} Locally holomorphic $p$-contact and $s$-symplectic manifolds are special cases of {\bf spin manifolds}. In particular, if $X$ is a compact complex manifold $X$ with $\mbox{dim}_\C X = 2p+1$ and $p$ is odd, the following implications hold: \begin{eqnarray}\label{Obs:p-contact_contact_spin_introd}X \hspace{2ex} \mbox{is {\bf holomorphic contact}}\hspace{2ex} \implies  X \hspace{2ex} \mbox{is {\bf holomorphic $p$-contact}}\hspace{2ex} \implies X \hspace{2ex} \mbox{is {\bf spin}}.\end{eqnarray}

\end{Obs}  

\vspace{1ex}

It was noted in [KPU25a, Observation 2.1.] that compact complex manifolds carrying a $\C$-valued holomorphic $p$-contact structure are never K\"ahler. However, locally holomorphic $p$-contact manifolds (i.e. those carrying a line-bundle-valued holomorphic $p$-contact structure) introduced in this paper as generalisations of their scalar counterparts can even be projective. 

After giving several basic properties of holomorphic $p$-contact manifolds in $\S$\ref{section:p-contact}, we investigate in $\S$\ref{section:pos-curvature} what kind of (possibly singular in some cases) Hermitian fibre metric the canonical bundle of a compact K\"ahler manifold $X$ of one of the three types in (\ref{Obs:p-contact_contact_spin_introd}) can possess. The hypotheses we make involve the relatively recent notion of partial positivity, called $m$-positivity of a $(1,\,1)$-current $T$ with respect to a given Hermitian metric $\omega$, introduced by Dieu in [Die06] and denoted in [DP25a], [DP25b] and in this paper by $T\geq_{m,\,\omega}0$. Theorems \ref{The:p-contact_n-p-pos}, \ref{The:contact_n-1-pos} and \ref{The:spin_scal-pos} are summarised in the following

\begin{The}\label{The:positivity_introd} Let $(X,\,\omega)$ be a compact {\bf K\"ahler} manifold with $\mbox{dim}_\C X = n$.

\vspace{1ex}

(i)\, Suppose there exists a {\bf spin structure} on $X$ and the anti-canonical bundle $-K_X$ possesses a $C^\infty$ Hermitian fibre metric $h^2$ such that \begin{eqnarray}\label{eqn:Scal_positivity_point}{\mbox Scal}\,(\omega,\,h^2)\geq 0 \hspace{2ex}\mbox{everywhere on}\hspace{1ex} X \hspace{3ex}\mbox{and}\hspace{3ex} {\mbox Scal}\,(\omega,\,h^2)> 0 \hspace{2ex}\mbox{at some point of}\hspace{1ex} X.\end{eqnarray} Then $H^0(X,\,K_X) = H^{n,\,0}(X,\,\C) = \{0\}$.

\vspace{1ex}

(ii)\, Suppose that $n = 2p+1$ and there exist a holomorphic line bundle $L$ on $X$ and an {\bf $L$-valued holomorphic contact structure} $\eta\in H^{1,\,0}_{\bar\partial}(X,\,L)$. Then, there exists {\bf no} (singular or otherwise) Hermitian fibre metric $h^{-(p+1)}$ on the canonical bundle $K_X\simeq L^{-(p+1)}$ of $X$ such that \begin{eqnarray*} i\Theta_{h^{-(p+1)}}(K_X)\geq_{n-1,\,\omega}0.\end{eqnarray*}

\vspace{1ex}

(iii)\, Suppose that $n = 2p+1$ and $p$ is odd. Further suppose there exist a holomorphic line bundle $L$ on $X$ and an {\bf $L$-valued holomorphic $p$-contact structure} $\Gamma\in H^{p,\,0}_{\bar\partial}(X,\,L)$. Then, there exists {\bf no} (singular or otherwise) Hermitian fibre metric $h^{-2}$ on the canonical bundle $K_X\simeq L^{-2}$ such that \begin{eqnarray*}i\Theta_{h^{-2}}(K_X)\geq_{n-p,\,\omega}0.\end{eqnarray*}

\end{The}

In part (i) above, ${\mbox Scal}\,(\omega,\,h^2)$ stands for the scalar curvature of $X$ with respect to $\omega$ and $h^2$ defined in the expected way in $\S$\ref{subsection:scalar-spin}. This result can be viewed as a version in our setting of a consequence (cf. e.g. [B\"ar11, Corollary 2.5.13.]) of Lichnerowicz's identity in [Lic63] (whose role is played here by the Bochner-Kodaira-Nakano identity) asserting that there are no non-trivial harmonic spinor fields on a connected compact Riemannian spin manifold with ${\mbox Scal}\geq 0$ everywhere and ${\mbox Scal}> 0$ somewhere.

Parts (ii) and (iii) above rule out $K_X$ being flat or positive but allow it to be negative. So, $X$ can very well be a Fano manifold. Fano contact manifolds were studied, for example, in [Bea98]. Fano $p$-contact manifolds are mentioned in (ii) of Corollary \ref{Cor:Demailly_K_not-psef} which translates into our language an observation of [Dem02].

An analogue of part (iii) above for locally holomorphic $s$-symplectic manifolds is given in Theorem \ref{The:s-symplectic_s-pos}: the curvature assumption is similar, but the conclusion is different. It merely asserts that the sheaf ${\cal F}_\Omega$ of germs of holomorphic $(1,\,0)$-vector fields $\xi$ that contract the bundle-valued holomorphic $s$-symplectic form $\Omega$ to 0 has rank zero.

These results, whose proofs use the local regularisation Theorem 1.6. in [DP25a] to handle singular fibre metrics, can be viewed as generalisations of the main theorem in [Dem02] that we restate as Theorem \ref{The:p-contact_n-p-pos}. Some of the consequences of this discussion are gathered in the following statement that merges parts of Corollaries \ref{Cor:Demailly_K_not-psef} and \ref{Cor:p-contact_non-Kobayashi-hyperbolic}.

\begin{Cor}\label{Cor:p-contact_consequences_introd} Let $X$ be a compact {\bf K\"ahler} manifold with $\mbox{dim}_\C X = n = 2p+1$ and $p$ odd. Suppose there exists an $L$-valued {\bf holomorphic $p$-contact structure} $\Gamma\in H^{p,\,0}_{\bar\partial}(X,\,L)$ for some holomorphic line bundle $L$ over $X$. Then:

\vspace{1ex}  

(i)\, the canonical line bundle $K_X$ of $X$ {\bf is not pseudo-effective}. Hence, the Kodaira dimension of $X$ is negative: $k(X) = -\infty$. 

\vspace{1ex}  

(ii)\, if $X$ is projective, $X$ is {\bf covered by rational curves}. In particular, $X$ is {\bf not Kobayashi hyperbolic}.  

\end{Cor}

In $\S$\ref{section:F-G_h}, we briefly discuss the versions of the sheaves ${\cal F}_\Gamma$ and ${\cal G}_\Gamma$ that played such a key role in the scalar setting of [KPU25a] and [KPU25b], updated now to our current bundle-valued structures where the latter sheaf depends on the Hermitian fibre metric $h$ on $L$.

In $\S$\ref{section:examples}, we discuss some examples of locally holomorphic $p$-contact manifolds. That discussion, combined with standard results, implies that the complex projective space $\Proj^n$ is: 

\vspace{1ex}

$\bullet$ never a (locally) holomorphic $s$-symplectic manifold;

\vspace{1ex}

$\bullet$ a spin manifold if and only if $n$ is odd (cf. e.g. [LM89, II, $\S2$, Example 2.4.]);

\vspace{1ex}

$\bullet$ a (locally) holomorphic $p$-contact manifold (for some $p$) if and only if $n\equiv 3$ mod $4$.

\vspace{1ex}

Moreover, Proposition \ref{Prop:hypersurfaces-projective-spaces} shows that the only smooth complex hypersurfaces in complex projective space that carry a holomorphic $p$-contact structure are copies of $\Proj^n$ with $n\equiv 3$ mod $4$, seen as hypersurfaces in $\Proj^{n+1}$.

Fortunately, Proposition \ref{Prop:p-contact_s-symplectic_fibration} shows that new examples of (locally) holomorphic $p$-contact manifolds can be obtained from existing ones via surjective holomorphic submersions $\pi:X\longrightarrow Y$ between compact complex manifolds with $\mbox{dim}_\C Y = 2s = 4r$ and $\mbox{dim}_\C X = 2s + 2p + 1$ when $p$ is odd. Loosely reworded, Proposition \ref{Prop:p-contact_s-symplectic_fibration} ensures that if the base manifold $Y$ carries a holomorphic {\bf $s$-symplectic structure} and if there exists a fibrewise holomorphic $p$-contact structure, then the total space $X$ of the fibration carries a holomorphic {\bf $(s+p)$-contact structure}. In particular, products $Y\times Z$ of a holomorphic $s$-symplectic manifold $Y$ and a holomorphic $p$-contact manifold $Z$ are holomorphic $(s+p)$-contact manifolds.

\vspace{2ex}

\noindent {\bf Acknowledgments}. Work on this project started during the second-named author's week-long visit to the University of Queensland in March 2025. He is grateful to the first-named author for the invitation and support. Thanks are also due to S. Dinew for valuable comments.

\section{The $p$-contact setting}\label{section:p-contact}


\subsection{Preliminaries}\label{subsection:p-contact_prelim}

Let $X$ be a compact complex manifold with $\mbox{dim}_\C X =n = 2p+1 = 4l+3$, where $p\geq 1$ and $l\geq 0$ are integers. Let $L$ be a holomorphic line bundle on $X$ and let $\Gamma\in C^\infty_{p,\,0}(X,\,L)$ be an $L$-valued $C^\infty$ $(p,\,0)$-form. We shall often assume $\Gamma$ to be holomorphic, namely $\bar\partial\Gamma = 0$ on $X$, where $\bar\partial$ is the canonical $(0,\,1)$-connection on $L$ induced by its holomorphic structure and extending to $L$-valued forms the Cauchy-Riemann operator $\bar\partial$ (denoted by the same symbol) acting on scalar-valued forms on $X$.

If we fix a collection of local holomorphic trivialisations $\theta_\alpha : L_{|U_\alpha}\longrightarrow U_\alpha\times\C$ of $L$ over an open covering ${\cal U} = (U_\alpha)_\alpha$ of $X$ and denote by $(e_\alpha)_\alpha$ the corresponding holomorphic frames over the $U_\alpha$'s and by $(g_{\alpha\beta})_{\alpha,\,\beta}$ the corresponding (invertible, holomorphic) transition functions $g_{\alpha\beta} : U_\alpha\cap U_\beta\longrightarrow\C$ of $L$, we have $e_\beta = g_{\alpha\beta}\,e_\alpha$ on $U_\alpha\cap U_\beta$ for all $\alpha,\beta$.

This implies that any Hermitian fibre metric $h$ on $L$ can be seen as a collection of locally defined functions $\varphi_\alpha : U_\alpha\longrightarrow\R$ satisfying the gluing condition \begin{eqnarray}\label{eqn:gluing_L-h}\varphi_\alpha = \varphi_\beta + \log|g_{\alpha\beta}|^2 \hspace{5ex} \mbox{on}\hspace{2ex} U_\alpha\cap U_\beta\end{eqnarray} for all $\alpha,\beta$. Indeed, the local weight functions $\varphi_\alpha$ of $h$ are defined by requiring \begin{eqnarray*}|e_\alpha|^2_h = e^{-\varphi_\alpha}\end{eqnarray*} on each $U_\alpha$. We will write $h=(e^{-\varphi_\alpha})_\alpha$. The metric $h$ is (standardly) said to be $C^\infty$ (resp. singular) if all the $\varphi_\alpha$'s are $C^\infty$ (resp. $L^1_{loc}$). When $h$ is $C^\infty$, the Chern connection $D_h$ of $(L,\,h)$ splits into a $(1,\,0)$-connection $D_h':C^\infty_{r,\,s}(X,\,L)\longrightarrow C^\infty_{r+1,\,s}(X,\,L)$ and a $(0,\,1)$-connection $D_h'' = \bar\partial:C^\infty_{r,\,s}(X,\,L)\longrightarrow C^\infty_{r,\,s+1}(X,\,L)$. So, $D_h = D_h' + \bar\partial$ and it is standard that\begin{eqnarray}\label{eqn:D'_h_U_alpha}(D'_h)_{|U_\alpha} = \partial - \partial\varphi_\alpha\wedge\cdot = e^{\varphi_\alpha}\,\partial\bigg(e^{-\varphi_\alpha}\,\cdot\,\bigg):=\partial_{-\varphi_\alpha} \hspace{5ex} \mbox{on}\hspace{2ex} U_\alpha\end{eqnarray} for every $\alpha$. It is equally standard that the resulting curvature form $i\Theta_h(L)\in C^\infty_{1,\,1}(X,\,\R)$, defined by $D'_h\bar\partial + \bar\partial D'_h = \Theta_h(L)\wedge\cdot\,$, is locally given by \begin{eqnarray*}i\Theta_h(L)_{|U_\alpha} = i\partial\bar\partial\varphi_\alpha\end{eqnarray*} for every $\alpha$. The last piece of standard information about Hermitian fibre metrics that we recall from [Dem97, chapter V] is the existence of a sesquilinear pairing \begin{eqnarray*}C^\infty_{r_1,\,s_1}(X,\,L)\times C^\infty_{r_2,\,s_2}(X,\,L)\longrightarrow C^\infty_{r_1+s_2,\,s_1+r_2}(X,\,\C)\end{eqnarray*} induced by $h$ on pairs $(\sigma,\,\tau)$ of global $L$-valued forms of any bidegrees. In a local trivialisation $(\theta_\alpha,\,e_\alpha)$ of $L$, it is defined by \begin{eqnarray}\label{eqn:h-pairing_def}\bigg\{\sigma_{|U_\alpha}=\sigma_\alpha\otimes e_\alpha,\,\tau_{|U_\alpha}=\tau_\alpha\otimes e_\alpha\bigg\}_h: = \sigma_\alpha\wedge\overline\tau_\alpha\,e^{-\varphi_\alpha}.\end{eqnarray} The compatibility of any connection $D$ on $L$ with a given metric $h$, as is the case with the Chern connection, is defined to mean that $D$ satisfies the Leibniz rule w.r.t. to this sesquilinear pairing whose result is a scalar form on $X$.

Now, for every $\Gamma\in C^\infty_{p,\,0}(X,\,L)$ and every $\alpha$, we have \begin{eqnarray*}\Gamma_{|U_\alpha} = \Gamma_\alpha\otimes e_\alpha \hspace{5ex} \mbox{on}\hspace{2ex} U_\alpha,\end{eqnarray*} where each $\Gamma_\alpha$ is a $C^\infty$ $(p,\,0)$-form on $U_\alpha$ and the local pieces satisfy the gluing condition \begin{eqnarray}\label{eqn:gluing_L-Gamma}\Gamma_\alpha = g_{\alpha\beta}\,\Gamma_\beta \hspace{5ex} \mbox{on}\hspace{2ex} U_\alpha\cap U_\beta\end{eqnarray} for all $\alpha,\beta$. Equivalently, $\Gamma$ can be seen as a collection $(\Gamma_\alpha)_\alpha$ of locally defined $C^\infty$ $\C$-valued $(p,\,0)$-forms satisfying (\ref{eqn:gluing_L-Gamma}). By taking $\bar\partial$ in (\ref{eqn:gluing_L-Gamma}) and using the holomorphicity of the $g_{\alpha\beta}$'s, we get $\bar\partial\Gamma_\alpha = g_{\alpha\beta}\,\bar\partial\Gamma_\beta$ on $U_\alpha\cap U_\beta$ for all  $\alpha,\beta$. Hence, $\bar\partial\Gamma = (\bar\partial\Gamma_\alpha)_\alpha\in C^\infty_{p,\,1}(X,\,L)$ is a global $L$-valued $C^\infty$ $(p,\,1)$-form on $X$. In particular, the condition $\bar\partial\Gamma = 0$ is equivalent to $\bar\partial\Gamma_\alpha = 0$ on $U_\alpha$ for every $\alpha$.

  However, $\partial\Gamma$ does not make global sense since there is no canonical $(1,\,0)$-connection on the holomorphic line bundle $L$. Another way of seeing this is to take $\partial$ in the equality (\ref{eqn:gluing_L-Gamma}) of scalar forms, which yields \begin{eqnarray*}\partial\Gamma_\alpha = g_{\alpha\beta}\,\partial\Gamma_\beta + \partial g_{\alpha\beta}\wedge\Gamma_\beta \hspace{5ex} \mbox{on}\hspace{2ex} U_\alpha\cap U_\beta\end{eqnarray*} for all  $\alpha,\beta$, and to observe that the last term prevents the locally defined scalar $(p+1,\,0)$-forms $\partial\Gamma_\alpha$ from gluing together into a global $L$-valued form. Nevertheless, $\Gamma_\beta\wedge\Gamma_\beta = 0$ since $p$ is assumed {\it odd}, hence we get \begin{eqnarray}\label{eqn:gluing_L-del-Gamma}\Gamma_\alpha\wedge\partial\Gamma_\alpha = g_{\alpha\beta}^2\,\Gamma_\beta\wedge\partial\Gamma_\beta \hspace{5ex} \mbox{on}\hspace{2ex} U_\alpha\cap U_\beta\end{eqnarray} for all $\alpha,\beta$. Since $(g_{\alpha\beta}^2)_{\alpha,\,\beta}$ are the transition functions of the holomorphic line bundle $L^2:=L\otimes L$, we infer the first statement in

\begin{Obs}\label{Obs:Gamma-wedge-sel-Gamma_global} If $\mbox{dim}_\C X = n= 2p+1$ and $p$ is {\bf odd}, for any $\Gamma\in C^\infty_{p,\,0}(X,\,L)$ we have:

\vspace{1ex}

(i)\, $\Gamma\wedge\partial\Gamma : = (\Gamma_\alpha\wedge\partial\Gamma_\alpha)_\alpha \in C^\infty_{n,\,0}(X,\,L^2)$ is a well-defined $C^\infty$ $L^2$-valued $(n,\,0)$-form on $X$;

\vspace{1ex}

(ii)\, if $h=(e^{-\varphi_\alpha})_\alpha$ is any $C^\infty$ Hermitian fibre metric on $L$ and $D_h = D_h' + \bar\partial$ is the associated Chern connection, we have \begin{eqnarray*}\Gamma\wedge\partial\Gamma = \Gamma\wedge D_h'\Gamma\in C^\infty_{n,\,0}(X,\,L^2).\end{eqnarray*}

\end{Obs}

\noindent {\it Proof.} To prove (ii), we first use (\ref{eqn:D'_h_U_alpha}) to get \begin{eqnarray*}(D'_h\Gamma)_{|U_\alpha} = (\partial\Gamma_\alpha - \partial\varphi_\alpha\wedge\Gamma_\alpha)\otimes e_\alpha\end{eqnarray*} for every $\alpha$. Hence, \begin{eqnarray*}(\Gamma\wedge D'_h\Gamma)_{|U_\alpha} = \bigg(\Gamma_\alpha\wedge(\partial\Gamma_\alpha - \partial\varphi_\alpha\wedge\Gamma_\alpha)\bigg)\otimes e_\alpha^2 = (\Gamma_\alpha\wedge\partial\Gamma_\alpha)\otimes e_\alpha^2\end{eqnarray*} for every $\alpha$ since $\Gamma_\alpha\wedge\Gamma_\alpha = 0$ due to $p$ being odd. This proves the contention. \hfill $\Box$

\subsection{Local holomorphic $p$-contact structures}\label{subsection:local_hol-p-contact}

The first main notion introduced in this paper is the following generalisation to an $L$-valued form of the notion of {\it holomorphic $p$-contact structure} introduced and studied for the case where $L$ is trivial in [KPU25]. In other words, we now introduce the local version of the global notion of [KPU25]. 
    
\begin{Def}\label{Def:L_hol-p-contact} Let $L$ be a holomorphic line bundle on a compact complex manifold $X$ with $\mbox{dim}_\C X = n = 2p+1 = 4l+3$, where $p$ and $l$ are integers.

\vspace{1ex}

$(1)$\, An {\bf $L$-valued holomorphic $p$-contact structure} on $X$ is a form $\Gamma\in C^\infty_{p,\,0}(X,\,L)$ satisfying the following two conditions: \begin{eqnarray}\label{eqn:L_hol-p-contact_Gamma}(a)\,\,\, \bar\partial\Gamma = 0  \hspace{5ex}\mbox{and}\hspace{5ex} (b)\,\,\,\Gamma\wedge\partial\Gamma \neq 0 \hspace{2ex}\mbox{at every point of}\hspace{2ex} X.\end{eqnarray}

When the line bundle $L$ is not specified, any such form $\Gamma$ is alternatively called a {\bf local holomorphic $p$-contact structure} on $X$.
  
\vspace{1ex}

$(2)$\, If there exists an $L$-valued holomorphic $p$-contact structure $\Gamma\in C^\infty_{p,\,0}(X,\,L)$, $L$ is called a {\bf holomorphic $p$-contact bundle}, $(X,\,L,\,\Gamma)$ is called a {\bf holomorphic $p$-contact triple} and $X$ is called a {\bf (locally holomorphic) $p$-contact manifold}.

\end{Def}

The initial observations are collected in the following

\begin{Prop}\label{Prop:initial_L_p-contact} Let $(X,\,L,\,\Gamma)$ be a {\bf holomorphic $p$-contact triple} with $\mbox{dim}_\C X = n = 2p+1 = 4l+3$. Then:

\vspace{1ex}

(i)\, $\Gamma$ induces a holomorphic line bundle {\bf isomorphism} $L^2\simeq -K_X$ between $L^2$ and the anti-canonical bundle of $X$. With respect to any family $(U_\alpha,\,\theta_\alpha,\,e_\alpha)_\alpha$ of local holomorphic trivialisations of $L$, the equivalent isomorphism $S_\Gamma : L^{-2}\longrightarrow K_X$ is defined by \begin{eqnarray*}e_\alpha^{-2}\longmapsto\Gamma_\alpha\wedge\partial\Gamma_\alpha \hspace{5ex}\mbox{on}\hspace{2ex} U_\alpha \end{eqnarray*} for every $\alpha$.

\vspace{1ex}

(ii)\, For every $C^\infty$ Hermitian fibre metric $h$ on $L$, let $h^2$ be the induced $C^\infty$ Hermitian fibre metric on $L^2$ and let $h^2_\Gamma$ be the $C^\infty$ Hermitian fibre metric on $-K_X$ induced by $h^2$ via the bundle isomorphism $L^2\simeq -K_X$ defined by $\Gamma$ under (i) above. Let $Ric(h^2):=\frac{i}{2\pi}\,\Theta_{h^2_\Gamma}(-K_X)$ be the curvature form of $-K_X$ w.r.t. $h^2_\Gamma$. If ${\cal H}_L$ is the set of Hermitian fibre metrics $h$ on $L$, the map \begin{eqnarray}\label{eqn:h-Ric_map}{\cal H}_L\ni h\longrightarrow Ric(h^2)\in c_1(X)_{BC}\end{eqnarray} is {\bf surjective}, where $c_1(X)_{BC}$ is the first Chern class of $X$ (namely the first Chern class of $T^{1,\,0}X$, or equivalently, the first Chern class of $-K_X$) in the sense of the Bott-Chern cohomology.  

Moreover, if $\sim$ is the equivalence relation on ${\cal H}_L$ defined by \begin{eqnarray*}h_1\sim h_2 \iff \exists\, C\in\R \hspace{2ex}\mbox{such that}\hspace{2ex} h_1 = e^{-C}\,h_2\end{eqnarray*} and $\widehat{h}$ is the equivalence class of any $h\in{\cal H}_L$, the map \begin{eqnarray}\label{eqn:h-Ric_map_bis}{\cal H}_L/\sim\ni \widehat{h}\longrightarrow Ric(h^2)\in c_1(X)_{BC}\end{eqnarray} is {\bf well defined} and {\bf bijective}.

\vspace{1ex}

(iii)\, For every $C^\infty$ Hermitian fibre metric $h$ on $L$, the scalar $(n,\,n)$-form \begin{eqnarray}\label{eqn:volume-form}dV_{\Gamma,\,h}:=i^{n^2}\,\bigg\{\Gamma\wedge\partial\Gamma,\,\Gamma\wedge\partial\Gamma\bigg\}_{h^2}\in C^\infty_{n,\,n}(X,\,\R)\end{eqnarray} is {\bf strictly positive} at every point of $X$. It will be called the {\bf volume form} induced by $\Gamma$ and $h$.

With respect to any family $(U_\alpha,\,\theta_\alpha,\,e_\alpha)_\alpha$ of local holomorphic trivialisations of $L$ in which $\Gamma = (\Gamma_\alpha)_\alpha$ and $h=(e^{-\varphi_\alpha})_\alpha$, we have \begin{eqnarray*}\Gamma_\alpha\wedge\overline\Gamma_\alpha\,e^{-\varphi_\alpha} = \Gamma_\beta\wedge\overline\Gamma_\beta\,e^{-\varphi_\beta}  \hspace{5ex}\mbox{on}\hspace{2ex} U_\alpha\cap U_\beta \hspace{2ex}\mbox{for all}\hspace{2ex} \alpha,\beta.\end{eqnarray*} Moreover, with respect to the $\C$-valued $(p,\,p)$-form $\Gamma\wedge\overline\Gamma\,e^{-\varphi}$ globally defined on $X$ by the requirement that its restriction to each $U_\alpha$ be $\Gamma_\alpha\wedge\overline\Gamma_\alpha\,e^{-\varphi_\alpha}$, the above volume form is given by \begin{eqnarray}\label{eqn:volume-form_formula}dV_{\Gamma,\,h}:=i\partial(\Gamma\wedge\overline\Gamma\,e^{-\varphi})\wedge\bar\partial(\Gamma\wedge\overline\Gamma\,e^{-\varphi}).\end{eqnarray}

\vspace{1ex}

(iv)\, Fix an arbitrary Hermitian metric $\omega$ on $X$. For every $C^\infty$ Hermitian fibre metric $h$ on $L$, let $C_h>0$ be the unique constant given by the generalisation in [TW10] of Yau's theorem of [Yau78] for which there exists a Hermitian metric $\omega_h = \omega + i\partial\bar\partial\varphi_h>0$ on $X$ such that \begin{eqnarray}\label{eqn:M-A_Yau_T-W}\omega_h^n = C_h\,dV_{\Gamma,\,h}.\end{eqnarray}

Such a metric $\omega_h$ is unique and satisfies \begin{eqnarray}\label{eqn:two-Riccis}Ric(h^2) = Ric(\omega_h),\end{eqnarray} where $Ric(\omega_h):=\frac{i}{2\pi}\Theta_{\omega_h}(-K_X)$ is the curvature form of $-K_X$ w.r.t. the Hermitian fibre metric induced  by $\omega_h^n$.

\end{Prop}

\noindent {\it Proof.} (i)\, Since $\Gamma\wedge\partial\Gamma$ is a non-vanishing holomorphic $(n,\,0)$-form with values in $L^2$, we have \begin{eqnarray*}\Gamma\wedge\partial\Gamma\in H^{n,\,0}_{\bar\partial}(X, L^2)\simeq H^0(X,\,\Omega^n\otimes L^2) = H^0(X,\,K_X\otimes L^2).\end{eqnarray*} Thus, $\Gamma\wedge\partial\Gamma$ is a non-vanishing global holomorphic section of the holomorphic line bundle $K_X\otimes L^2$. The existence of such a section implies the triviality of $K_X\otimes L^2$, hence the isomorphism $L^2\simeq -K_X$.

Moreover, the non-vanishing section $\Gamma\wedge\partial\Gamma$ of $K_X\otimes L^2$ identifies in any local trivialisation the local holomorphic frame $e_\alpha^{-2}$ of $L^{-2}$ on $U_\alpha$ with the local holomorphic frame $\Gamma_\alpha\wedge\partial\Gamma_\alpha$ of $K_X$.

\vspace{1ex}

(ii)\, Let $\widetilde\omega$ be an arbitrary $C^\infty$ $(1,\,1)$-form lying in the class $c_1(X)_{BC} = c_1(-K_X)_{BC}$. There exists a $C^\infty$ Hermitian fibre metric $\widetilde{h}$ on $-K_X$ whose curvature form is $\widetilde\omega$ (up to a multiplicative constant): \begin{eqnarray*}\frac{i}{2\pi}\,\Theta_{\widetilde{h}}(-K_X) = \widetilde\omega.\end{eqnarray*} Since $-K_X\simeq L^2$, $\widetilde{h}$ induces a $C^\infty$ Hermitian fibre metric, that we still denote by $\widetilde{h}$, on $L^2$. Let $(\psi_\alpha)_\alpha$ be the local weight functions defining $\widetilde{h}$ in a collection of local trivialisations for $L$, namely $\widetilde{h} = (e^{-\psi_\alpha})_\alpha$. Then $\widetilde{h}^{\frac{1}{2}}:= (e^{-\frac{1}{2}\,\psi_\alpha})_\alpha$ is a $C^\infty$ Hermitian fibre metric on $L$. By construction, we have \begin{eqnarray*}Ric\bigg(\bigg(\widetilde{h}^{\frac{1}{2}}\bigg)^2\bigg) = Ric\bigg(\widetilde{h}\bigg) = \widetilde\omega.\end{eqnarray*} This proves the surjectivity of the map (\ref{eqn:h-Ric_map}).

Now, fix an arbitrary $C^\infty$ Hermitian fibre metric $h_0=(e^{-\psi_\alpha})_\alpha$ on $L$. Then $\psi_\alpha = \psi_\beta + \log|g_{\alpha\beta}|^2$ on $U_\alpha\cap U_\beta$ for all $\alpha, \beta$. If $h=(e^{-\varphi_\alpha})_\alpha$ is any $C^\infty$ Hermitian fibre metric on $L$, we have $\varphi_\alpha = \varphi_\beta + \log|g_{\alpha\beta}|^2$ on $U_\alpha\cap U_\beta$ for all $\alpha, \beta$. Hence, \begin{eqnarray*}\varphi_\alpha - \psi_\alpha = \varphi_\beta - \psi_\beta  \hspace{5ex}\mbox{on}\hspace{2ex} U_\alpha\cap U_\beta\end{eqnarray*} for all $\alpha,\beta$. Therefore, the exists a global $C^\infty$ function $f:X\longrightarrow\R$ whose restriction to each $U_\alpha$ is $\varphi_\alpha - \psi_\alpha$. We conclude that every $C^\infty$ Hermitian fibre metric $h$ on $L$ can be specified in terms of a reference such metric $h_0=(e^{-\psi_\alpha})_\alpha$ by means of such a function $f$: \begin{eqnarray}\label{eqn:h-h_0_f}h = (e^{-f-\psi_\alpha})_\alpha = e^{-f}\,h_0.\end{eqnarray}

  Moreover, for any two $C^\infty$ Hermitian fibre metrics $h_1= e^{-f_1}\,h_0$ and $h_2= e^{-f_2}\,h_0$ on $L$, we have the following equivalences (in which we identify $h_l^2$ with $h_{l,\,\Gamma}^2$ for $l=0, 1,2$) for equalities on $X$ : \begin{eqnarray*}Ric(h_1^2) = Ric(h_2^2) & \iff & \frac{i}{2\pi}\,\Theta_{h_1^2}(-K_X) = \frac{i}{2\pi}\,\Theta_{h_2^2}(-K_X) \\
    & \iff & \frac{i}{2\pi}\,\Theta_{h_0^2}(-K_X) + 2\,i\partial\bar\partial f_1 = \frac{i}{2\pi}\,\Theta_{h_0^2} + 2\,i\partial\bar\partial f_2 \\
  & \iff & \partial\bar\partial(f_1 - f_2) = 0  \iff f_1 - f_2 = \mbox{Const}.\end{eqnarray*} We have used the compactness of $X$ to get the constancy conclusion for $f_1 - f_2$. This proves that $Ric(h_1^2) = Ric(h_2^2)$ if and only if there exists a constant $C\in\R$ such that $h_1 = \exp(-C)\,h_2$, which proves the well-definedness (i.e. its independence of the choice of representative of the equivalence class $\widehat{h}$) and the injectivity of the map (\ref{eqn:h-Ric_map_bis}).

\vspace{1ex}

(iii)\, Since $\Gamma\wedge\partial\Gamma = (\Gamma_\alpha\wedge\partial\Gamma_\alpha)_\alpha$ is an $L^2$-valued $(n,\,0)$-form on $X$ (see Observation \ref{Obs:Gamma-wedge-sel-Gamma_global}) and since $h^2 = (e^{-2\varphi_\alpha})_\alpha$ is a fibre metric on $L^2$ whenever $h = (e^{-\varphi_\alpha})_\alpha$ is one on $L$, the general definition (\ref{eqn:h-pairing_def}) yields: \begin{eqnarray}\label{eqn:volume-form_local-expression}(dV_{\Gamma,\,h})_{|U_\alpha} = i^{n^2}\,(\Gamma_\alpha\wedge\partial\Gamma_\alpha)\wedge(\overline\Gamma_\alpha\wedge\partial\overline\Gamma_\alpha)\,e^{-2\varphi_\alpha}.\end{eqnarray} Now, it is standard (see e.g. [Dem97, III.$\S1.$A]) that $i^{n^2}\,(\Gamma_\alpha\wedge\partial\Gamma_\alpha)\wedge(\overline\Gamma_\alpha\wedge\partial\overline\Gamma_\alpha)\geq 0$ as an $(n,\,n)$-form for every $(n,\,0)$-form $\Gamma_\alpha$. In our case, $\Gamma_\alpha\wedge\partial\Gamma_\alpha$ is non-vanishing, so the above inequality is strict.

As for the next claim, putting together the gluing properties (\ref{eqn:gluing_L-Gamma}) and (\ref{eqn:gluing_L-h}), we get: \begin{eqnarray*}\Gamma_\alpha\wedge\overline\Gamma_\alpha\,e^{-\varphi_\alpha} = \bigg(|g_{\alpha\beta}|^2\,\Gamma_\beta\wedge\overline\Gamma_\beta\bigg)\,\bigg(e^{-\varphi_\beta}\,\frac{1}{|g_{\alpha\beta}|^2}\bigg)\hspace{5ex}\mbox{on}\hspace{2ex} U_\alpha\cap U_\beta \hspace{2ex}\mbox{for all}\hspace{2ex} \alpha,\beta,\end{eqnarray*} as claimed. Moreover, since $\partial\overline\Gamma_\alpha = 0$ and $\bar\partial\Gamma_\alpha = 0$, we get: \begin{eqnarray*}\partial(\Gamma_\alpha\wedge\overline\Gamma_\alpha\,e^{-\varphi_\alpha}) & = & \partial\Gamma_\alpha\wedge\overline\Gamma_\alpha\,e^{-\varphi_\alpha} - \Gamma_\alpha\wedge\overline\Gamma_\alpha\,(e^{-\varphi_\alpha}\,\partial\varphi_\alpha) \\
  \bar\partial(\Gamma_\alpha\wedge\overline\Gamma_\alpha\,e^{-\varphi_\alpha}) & = & -\Gamma_\alpha\wedge\bar\partial\overline\Gamma_\alpha\,e^{-\varphi_\alpha} - \Gamma_\alpha\wedge\overline\Gamma_\alpha\,(e^{-\varphi_\alpha}\,\bar\partial\varphi_\alpha).\end{eqnarray*} Multiplying these two equalities and using the fact that $\Gamma_\alpha\wedge\Gamma_\alpha = 0$ and $\overline\Gamma_\alpha\wedge\overline\Gamma_\alpha = 0$ (since $p$ is odd), we get: \begin{eqnarray*}i\partial(\Gamma_\alpha\wedge\overline\Gamma_\alpha\,e^{-\varphi_\alpha})\wedge\bar\partial(\Gamma_\alpha\wedge\overline\Gamma_\alpha\,e^{-\varphi_\alpha}) & = & -i\,\partial\Gamma_\alpha\wedge\overline\Gamma_\alpha\wedge\Gamma_\alpha\wedge\bar\partial\overline\Gamma_\alpha\,e^{-2\varphi_\alpha} \\
  & = & i\,(\Gamma_\alpha\wedge\partial\Gamma_\alpha)\wedge(\overline\Gamma_\alpha\wedge\partial\overline\Gamma_\alpha)\,e^{-2\varphi_\alpha} = (dV_{\Gamma,\,h})_{|U_\alpha},\end{eqnarray*} where the last equality is (\ref{eqn:volume-form_local-expression}) since $i^{n^2} = i$ due $n$ being odd. This proves (\ref{eqn:volume-form_formula}).

\vspace{1ex}

(iv)\, The uniqueness of a solution $\omega_h = \omega + i\partial\bar\partial\varphi_h>0$ of the Monge-Amp\`ere equation (\ref{eqn:M-A_Yau_T-W}) is guaranteed by the Tosatti-Weinkove generalisation in [TW10] to the Hermitian case of Yau's celebrated theorem for K\"ahler metrics of [Yau78]. It remains to prove equality (\ref{eqn:two-Riccis}).

Let $h = (e^{-\varphi_\alpha})_\alpha$ an arbitrary $C^\infty$ Hermitian fibre metric on $L$. The metric $h^2 = (e^{-2\varphi_\alpha})_\alpha$ it induces on $L^2$ induces, in turn, a fibre metric $h_\Gamma^2$ on $-K_X$ that makes the isomorphism $L^2\simeq -K_X$ into an isometry. The associated curvature form is given in any local trivialisation by \begin{eqnarray}\label{eqn:Ric_h2_proof_equality_1}Ric(h^2)_{|U_\alpha} = \frac{i}{2\pi}\,\Theta_{h^2_\Gamma}(-K_X)_{|U_\alpha} = \frac{i}{2\pi}\,\partial\bar\partial(2\varphi_\alpha) = \frac{i}{\pi}\,\partial\bar\partial\varphi_\alpha\end{eqnarray} for all $\alpha$.

  On the other hand, the isometry $S_\Gamma:(L^{-2},\,h^{-2} = (e^{2\varphi_\alpha})_\alpha)\simeq (K_X,\, h_\Gamma^{-2})$, thanks to its shape given under (i) above, yields the latter equality in \begin{eqnarray*}e^{\varphi_\alpha} = |e_\alpha^{-2}|_{h^{-2}} = |\Gamma_\alpha\wedge\partial\Gamma_\alpha|_{h_\Gamma^{-2}},\end{eqnarray*} hence $\varphi_\alpha = \log|\Gamma_\alpha\wedge\partial\Gamma_\alpha|_{h_\Gamma^{-2}}$ on $U_\alpha$ for every $\alpha$. We get: \begin{eqnarray}\label{eqn:Ric_h2_proof_equality_2}\nonumber Ric(\omega_h)_{|U_\alpha} & = & -\frac{i}{2\pi}\,\partial\bar\partial\log(\omega_h^n)_{|U_\alpha} = -\frac{i}{2\pi}\,\partial\bar\partial\log\bigg(C_h\,\frac{dV_{\Gamma,\,h}}{dV_n}\bigg)_{|U_\alpha} = \frac{i}{2\pi}\,\partial\bar\partial\log|\Gamma_\alpha\wedge\partial\Gamma_\alpha|^2_{h_\Gamma^{-2}} \\
    & = & \frac{i}{2\pi}\,\partial\bar\partial(2\varphi_\alpha) = \frac{i}{\pi}\,\partial\bar\partial\varphi_\alpha,\end{eqnarray} where we used the local shape (\ref{eqn:volume-form_local-expression}) on $U_\alpha$ of $dV_{\Gamma,\,h}$ to get the third equality in which $dV_n:=idz_1\wedge d\bar{z}_1\wedge\dots\wedge idz_n\wedge d\bar{z}_n$ is the standard volume form defined by local coordinates on $U_\alpha$.

  The contention follows by putting together (\ref{eqn:Ric_h2_proof_equality_1} ) and (\ref{eqn:Ric_h2_proof_equality_2}). \hfill $\Box$

\subsection{Links with contact and spin manifolds}\label{subsection:links_contact-spin} 

The main takeaway from (i) of Proposition \ref{Prop:initial_L_p-contact} is that the existence of a holomorphic $p$-contact structure on $X$ forces the anti-canonical bundle (hence also the canonical bundle) of $X$ to have a square root (in the person of the holomorphic $p$-contact bundle $L$, respectively its dual $L^{-1}$). In other words, we get

\begin{Cor}\label{Cor:p-contact_spin} Any compact complex manifold $X$ with $\mbox{dim}_\C X \equiv 3 \hspace{2ex} mod \hspace{1ex} 4$ that carries a {\bf local holomorphic $p$-contact structure} is a {\bf spin manifold}.

\end{Cor}

Recall (see e.g. the book [B\"ar11] for details) that a spin manifold is an oriented Riemannian manifold $(M,\,g)$ equipped with a {\it spin structure} $(P^{Spin}(M),\,\bar\rho)$, namely a pair consisting of a {\it principal $Spin(m)$-bundle} $P^{Spin}(M)$ over $M$ (where $m = \dim_\R M$) and a {\it double covering} $\bar\rho:P^{Spin}(M)\longrightarrow P^{SO}(M)$, where $P^{SO}(M)$ is the principal $SO(m)$-bundle of oriented orthonormal frames of the tangent bundle $TM$ of $M$, such that the restriction of $\bar\rho$ to each fibre of $P^{Spin}(M)$ is the standard double covering $\rho:Spin(m)\longrightarrow SO(m)$. Isomorphic double coverings are considered to define a same spin structure.

Recall that the spin group $Spin(m)\subset Cl_m$ is defined as the subset consisting of Clifford products of an even number of points on the unit sphere $S^{m-1}\subset\R^m$ of the Clifford algebra $Cl_m$ for $\R^m$ equipped with the Euclidean inner product. The set $Spin(m)$ is a group with respect to the Clifford multiplication of $Cl_m$. It is even a Lie group that is compact for all $m\in\N$, connected for all $m\geq 2$ and simply connected for all $m\geq 3$. It is equipped with the group morphism $\rho:Spin(m)\longrightarrow SO(m)$ defined, for every $a\in Spin(m)$, as the element $\rho(a):\R^m\longrightarrow\R^m$ of $SO(m)$ given by \begin{eqnarray*}\rho(a) x:=a\cdot x\cdot a^{-1}, \hspace{5ex} x\in\R^n,\end{eqnarray*} where $\cdot$ denotes Clifford multiplication in $\R^m\subset Cl_m$. The map $\rho:Spin(m)\longrightarrow SO(m)$ is surjective and every element in $SO(m)$ has exactly two preimages in $Spin(m)$. That these two preimages can be thought of as the two square roots of a same element in $SO(m)$ is made explicit by the fact that, when $m=2$, one has $Spin(2)\simeq SO(2)\simeq S^1$, where $S^1\subset\C$ is the unit circle, and the map $\rho:Spin(2)\longrightarrow SO(2)$ is given by the squaring map $S^1\ni z\mapsto z^2\in S^1$.

Thus, a spin structure $(P^{Spin}(M),\,\bar\rho)$ on $(M,\,g)$, when it exists, can be thought of as a square root bundle of the principal $SO(m)$-bundle $P^{SO}(M)$. When the $m$-dimensional Riemannian manifold $(M,\,g)$ is replaced by a compact complex Hermitian manifold $(X,\,\omega)$ with $\mbox{dim}_\C X = n$, there is a nice criterion equating the existence of a spin structure on $X$ to the existence of a square root of the canonical bundle of $X$ (or, equivalently, of the anti-canonical  bundle). 

\begin{Prop}([Ati71, Proposition 3.2], [LM89])\label{Prop:spin_existence_standard} Let $X$ be a compact complex manifold and let $K_X$ be its canonical bundle.

\vspace{1ex}

(1)\, There exists a {\bf spin structure} on $X$ if and only if there exists a holomorphic line bundle $E$ on $X$ such that $E^2\simeq K_X$.

\vspace{1ex}

(2)\, The {\bf spin structures} on $X$ correspond bijectively to the isomorphism classes of holomorphic line bundles $E$ on $X$ such that $E^2\simeq K_X$.

\end{Prop}

\vspace{3ex} 

Another observation is that the notion of (local) holomorphic $p$-contact structure generalises the standard notion of {\it (local) holomorphic contact structure}. Recall that if $F$ is a holomorphic line bundle on a complex manifold $X$ with $\mbox{dim}_\C X = n = 2p+1$, an {\it $F$-valued holomorphic contact structure} on $X$ is a form $\eta\in C^\infty_{1,\,0}(X,\,F)$ satisfying the following two conditions: \begin{eqnarray}\label{eqn:L_hol-contact_eta}(a)\,\,\, \bar\partial\eta = 0  \hspace{5ex}\mbox{and}\hspace{5ex} (b)\,\,\,\eta\wedge(\partial\eta)^p \neq 0 \hspace{2ex}\mbox{at every point of}\hspace{2ex} X.\end{eqnarray}

As we saw earlier with $\Gamma\in C^\infty_{p,\,0}(X,\,L)$, one sees at once that $\partial\eta$ does not make sense as a vector-bundle-valued $(2,\,0)$-form on $X$, but the exterior product $\eta\wedge(\partial\eta)^p\in H^{n,\,0}(X,\,F^{p+1})$ is a well-defined holomorphic $F^{p+1}$-valued $(n,\,0)$-form on $X$. This amounts to $\eta\wedge(\partial\eta)^p$ being a global holomorphic section of the holomorphic line bundle $K_X\otimes F^{p+1}$. Since the {\it contact} definition (\ref{eqn:L_hol-contact_eta}) also imposes it to be non-vanishing, the homorphic line bundle $K_X\otimes F^{p+1}$ must be trivial. Equivalently, $F$ is a $(p+1)^{st}$ root of the anti-canonical bundle of $X$: \begin{eqnarray}\label{eqn:contact_root_p+1}F^{p+1}\simeq -K_X.\end{eqnarray}

\begin{Prop}\label{Prop:1-contact_p-contact} Let $F$ be a holomorphic line bundle on a compact complex manifold $X$ with $\mbox{dim}_\C X = n = 2p+1 = 4l+3$, where $p$ and $l$ are integers.

  If $\eta\in C^\infty_{1,\,0}(X,\,F)$ is an $F$-valued holomorphic {\bf contact structure} on $X$, then $\Gamma:=\eta\wedge(\partial\eta)^l\in H^{p,\,0}_{\bar\partial}(X,\,F^{l+1})$ is a well-defined $L:= F^{l+1}$-valued holomorphic {\bf $p$-contact structure} on $X$.

\end{Prop}

\noindent {\it Proof.} Let $((U_\alpha)_\alpha,\,(h_{\alpha\beta})_{\alpha,\,\beta})$ be a collection of (holomorphic) transition functions $h_{\alpha\beta}:U_\alpha\cap U_\beta\longrightarrow\C$ associated with a family of local holomorphic trivialisations of $F$ over open subsets $U_\alpha\subset X$. With respect to it, the $F$-valued holomorphic contact form can be viewed as a family $\eta=(\eta_\alpha)_\alpha$ of $\C$-valued locally defined holomorphic $(1,\,0)$-forms $\eta_\alpha\in H^{1,\,0}_{\bar\partial}(U_\alpha,\,\C)$ connected to one another by the property: \begin{eqnarray*}\eta_\alpha = h_{\alpha\beta}\,\eta_\beta  \hspace{6ex}\mbox{on}\hspace{2ex} U_\alpha\cap U_\beta,\end{eqnarray*} for all $\alpha,\beta$. Taking $\partial$, we get: \begin{eqnarray*}\partial\eta_\alpha & = & h_{\alpha\beta}\,\partial\eta_\beta + \partial h_{\alpha\beta}\wedge\eta_\beta, \\
  \eta_\alpha\wedge\partial\eta_\alpha & = & h_{\alpha\beta}^2\,\eta_\beta\wedge\partial\eta_\beta    \hspace{6ex}\mbox{on}\hspace{2ex} U_\alpha\cap U_\beta,\end{eqnarray*} for all $\alpha,\beta$, where the last equality follows from $\eta_\beta\wedge\eta_\beta = 0$ due to $\eta_\beta$ being a $1$-form. By induction on the power of $\partial\eta_\alpha$, we then get: \begin{eqnarray*}\eta_\alpha\wedge(\partial\eta_\alpha)^l & = & h_{\alpha\beta}^{l+1}\,\eta_\beta\wedge(\partial\eta_\beta)^l    \hspace{6ex}\mbox{on}\hspace{2ex} U_\alpha\cap U_\beta,\end{eqnarray*} for all $\alpha,\beta$. This shows that the locally defined holomorphic $\C$-valued $(p,\,0)$-forms $\Gamma_\alpha:=\eta_\alpha\wedge(\partial\eta_\alpha)^l$ glue together into an $F^{l+1}$-valued global holomorphic $(p,\,0)$-form $\Gamma:=(\Gamma_\alpha)_\alpha$ that we denote by $\eta\wedge(\partial\eta)^l$.

We get: \begin{eqnarray*}\Gamma\wedge\partial\Gamma & = & (\Gamma_\alpha\wedge\partial\Gamma_\alpha)_\alpha = \bigg(\eta_\alpha\wedge(\partial\eta_\alpha)^p\bigg)_\alpha = \eta\wedge(\partial\eta)^p\end{eqnarray*} and this form is non-vanishing thanks to $\eta$ being a holomorphic contact structure on $X$ (see (b) of (\ref{eqn:L_hol-contact_eta})). This proves the contention.   \hfill $\Box$

  \vspace{3ex}

  The upshot of these observations is that the class of holomorphic $p$-contact manifolds is intermediate between the class of holomorphic contact manifolds and that of spin manifolds (cf. second part of Observation \ref{Obs:p-contact_s-symplectic_spin_introd}).



 \subsection{Space of all $p$-contact structures}\label{subsection:all-p-contact} Let $X$ be a compact complex manifold with $\mbox{dim}_\C X = 2p+1$ and $p$ odd. Suppose $X$ is a spin manifold and fix a spin structure thereon, identified with a holomorphic line bundle $L$ on $X$ such that $L^2\simeq -K_X$. Then, \begin{eqnarray*}H^{n,\,0}_{\bar\partial}(X,\,L^2)\simeq H^0(X,\,K_X\otimes L^2)\simeq\C,\end{eqnarray*} the last isomorphism being a consequence of the triviality of $K_X\otimes L^2$ and the compactness of $X$.

We get the following holomorphic map \begin{eqnarray}\label{eqn:gamma-map-homogeneous}H^{p,\,0}_{\bar\partial}(X,\,L)\ni\Gamma\stackrel{T}{\longmapsto}\Gamma\wedge\partial\Gamma\in\C.\end{eqnarray} It is homogeneous of degree $2$ on $\C^{N+1}\simeq H^{p,\,0}_{\bar\partial}(X,\,L)$, where $N+1:=\dim_\C H^{p,\,0}_{\bar\partial}(X,\,L)$ (supposing that $\dim_\C H^{p,\,0}_{\bar\partial}(X,\,L)\geq 1$), hence it must be a homogeneous polynomial of degree $2$ thanks to the Cauchy inequalities being satisfied by holomorphic maps and to its domain of definition being the whole of a $\C^{N+1}$. It follows that $T$ identifies with a global holomorphic section of the holomorphic line bundle ${\cal O}(2)$ over the complex projective space $\Proj(H^{p,\,0}_{\bar\partial}(X,\,L))\simeq\Proj^N$ of $\C$-lines in $H^{p,\,0}_{\bar\partial}(X,\,L)$. Meanwhile, a form $\Gamma\in H^{p,\,0}_{\bar\partial}(X,\,L)$ is a holomorphic $p$-contact structure on $X$ if and only if $T(\Gamma)\neq 0$.

\begin{Conc}\label{Conc:all-p-contact_L} Let $L$ be a holomorphic line bundle over a compact complex manifold $X$ with $\mbox{dim}_\C X = 2p+1$ and $p$ odd such that $L^2\simeq -K_X$.

  Then, the set of $L$-valued holomorphic $p$-contact structures $\Gamma$ on $X$ is either empty or the complement of a degree-$2$ complex hypersurface in the complex projective space $\Proj(H^{p,\,0}_{\bar\partial}(X,\,L))$.

\end{Conc}

\subsection{Local holomorphic $p$-no-contact structures}\label{subsection:local_hol-p-no-contact}

As in the scalar case of [KPU25], we will sometimes contrast our local holomorphic $p$-contact structures with the opposite notion introduced in

\begin{Def}\label{Def:local_hol-p-no-contact} Let $L$ be a holomorphic line bundle on a compact complex manifold $X$ with $\mbox{dim}_\C X = n = 2p+1 = 4l+3$, where $p$ and $l$ are integers. Let $h$ be a $C^\infty$ Hermitian fibre metric on $L$ and let $D_h=D_h' + \bar\partial$ be the Chern connection of $(L,\,h)$.

\vspace{1ex}

An {\bf $(L,\,h)$-valued holomorphic $p$-no-contact structure} on $X$ is a form $\Gamma\in C^\infty_{p,\,0}(X,\,L)$ satisfying the following two conditions: \begin{eqnarray}\label{eqn:L_hol-p-no-contact_Gamma}(a)\,\,\, \bar\partial\Gamma = 0 \hspace{5ex}\mbox{and}\hspace{5ex} (b)\,\,\, D_h'\Gamma = 0.\end{eqnarray}

When the Hermitian line bundle $(L,\,h)$ is not specified, any such form $\Gamma$ is alternatively called a {\bf local holomorphic $p$-no-contact structure} on $X$.

\end{Def}  

The $(L,\,h)$-valued $p$-no-contact condition (\ref{eqn:L_hol-p-no-contact_Gamma}) on $\Gamma$ is equivalent to the condition $D_h\Gamma = 0$, i.e. to $\Gamma$ being {\bf $D_h$-parallel}, which amounts to (see (\ref{eqn:D'_h_U_alpha})) \begin{eqnarray}\label{eqn:L_hol-p-no-contact_Gamma_bis}\bar\partial\Gamma_\alpha = 0 \hspace{5ex}\mbox{and}\hspace{5ex} \partial\Gamma_\alpha = \partial\varphi_\alpha\wedge\Gamma_\alpha  \hspace{5ex}\mbox{on}\hspace{2ex} U_\alpha \hspace{2ex}\mbox{for every}\hspace{2ex} \alpha.\end{eqnarray} This implies that $\Gamma_\alpha\wedge\partial\Gamma_\alpha = -\partial\varphi_\alpha\wedge(\Gamma_\alpha\wedge\Gamma_\alpha) = 0$ since $p$ is odd, so $(L,\,h)$-valued holomorphic $p$-no-contact condition is incompatible with the $L$-valued holomorphic $p$-contact condition.

\section{Positivity and curvature}\label{section:pos-curvature} In this section, we study the types of smooth and singular fibre metrics that the canonical bundles of the manifolds discussed in this paper can carry.

\subsection{Laplacians induced by singular fibre metrics}\label{subsection:Laplacians_singular} Let $L'$ be a holomorphic line bundle on an arbitrary complex manifold $X$. It is standard that every (possibly singular) Hermitian fibre metric $\tilde{h}=(e^{-\psi_\alpha})_\alpha$ on $L'$ can be specified in terms of an arbitrarily fixed such metric $h=(e^{-\varphi_\alpha})_\alpha$ (that can be supposed $C^\infty$) as $\tilde{h} = he^{-f}$, where $f:X\longrightarrow\R\cup\{-\infty\}$ is a globally defined upper semicontinuous $L^1_{loc}$ function. This follows at once from the gluing condition (\ref{eqn:gluing_L-h}) that is satisfied by the two sets $(\varphi_\alpha)_\alpha$ and $(\psi_\alpha)_\alpha$ of local weight functions. Indeed, it implies that $\psi_\alpha - \varphi_\alpha = \psi_\beta - \varphi_\beta$ on $U_\alpha\cap U_\beta$ for all $\alpha$ and $\beta$, which proves the existence of a global function $f$ on $X$ such that $f = \psi_\alpha - \varphi_\alpha$ on $U_\alpha$ for every $\alpha$.      

    Fix an arbitrary Hermitian metric $\omega$ on $X$. Standard computations (that are left to the reader) prove explicit formulae for the Chern connection $D_{he^{-f}} = D'_{he^{-f}} + \bar\partial$ of $(L',\,{he^{-f}})$, the formal adjoints $(D'_{he^{-f}})^\star_{\omega,\,he^{-f}}$ and $\bar\partial^\star_{\omega,\,he^{-f}}$ of $D'_{he^{-f}}$, respectively $\bar\partial$, with respect to the $L^2$-inner product induced by $\omega$ and $he^{-f}$, as well as for the induced Laplacians: \begin{eqnarray*}\Delta'_{\omega,\,he^{-f}} & := & D'_{he^{-f}}\,(D'_{he^{-f}})^\star_{\omega,\,he^{-f}} + (D'_{he^{-f}})^\star_{\omega,\,he^{-f}}D'_{he^{-f}}:C^\infty_{p,\,q}(X,\,L')\longrightarrow C^\infty_{p,\,q}(X,\,L') \\
      \Delta''_{\omega,\,he^{-f}} & := & \bar\partial\,\bar\partial^\star_{\omega,\,he^{-f}} + \bar\partial^\star_{\omega,\,he^{-f}}\bar\partial:C^\infty_{p,\,q}(X,\,L')\longrightarrow C^\infty_{p,\,q}(X,\,L')\end{eqnarray*} when $f$ is $C^\infty$. We adopt those formulae as definitions in the general case of a possibly singular $f$, thereby obtaining differential operators on spaces of $L'$-valued $C^\infty$ forms with values in spaces of $L'$-valued currents. Specifically, we set:

    \begin{Def}\label{Def:Laplacians_singular} Let $(L',\,h)\longrightarrow(X,\,\omega)$ be a holomorphic line bundle equipped with a $C^\infty$ Hermitian fibre metric on a complex Hermitian manifold with $\mbox{dim}_\C X = n$. Let $f:X\longrightarrow\R\cup\{-\infty\}$ be an upper semicontinuous $L^1_{loc}$ function.

      For all $p,q\in\{0,\dots , n\}$, we define the following differential operators with values in vector spaces ${\cal D}'_{\,\bullet\,,\,\bullet\,}(X,\,L)$ of $L'$-valued currents of the specified bidegrees on $X$: \begin{eqnarray}\label{eqn:Laplacians_singular}\nonumber D'_{he^{-f}}:C^\infty_{p,\,q}(X,\,L') \longrightarrow  {\cal D}'_{p+1,\,q}(X,\,L'),  \hspace{3ex} D'_{he^{-f}} = \partial - \partial(\varphi_\alpha + f)\wedge\cdot = e^{\varphi_\alpha + f}\,\partial\bigg(e^{-\varphi_\alpha - f}\cdot\,\bigg) \end{eqnarray} \vspace{-5ex} \begin{eqnarray}\nonumber\bigg(D'_{he^{-f}}\bigg)^\star_{\omega,\,he^{-f}} & = & \partial^\star_\omega:C^\infty_{p,\,q}(X,\,L')  \longrightarrow  C^\infty_{p-1,\,q}(X,\,L')\\
        D''_{he^{-f}} & = & \bar\partial:C^\infty_{p,\,q}(X,\,L')  \longrightarrow  C^\infty_{p,\,q+1}(X,\,L')\\
 \nonumber       \bar\partial^\star_{\omega,\,he^{-f}} =  \bigg(\bar\partial + \bar\partial(\varphi_\alpha + f)\wedge\cdot\,\bigg)^\star_\omega & = & e^{\varphi_\alpha + f}\,\bar\partial^\star_\omega\bigg(e^{-\varphi_\alpha - f}\cdot\,\bigg):C^\infty_{p,\,q}(X,\,L')  \longrightarrow  {\cal D}'_{p,\,q-1}(X,\,L')\\
 \nonumber        \Delta'_{\omega,\,he^{-f}} & = & \Delta'_{\omega,\,h} - \bigg[[\Lambda_\omega,\,\bar\partial],\,i\partial f\wedge\cdot\,\bigg]:C^\infty_{p,\,q}(X,\,L')\longrightarrow {\cal D}'_{p,\,q}(X,\,L') \\
 \nonumber       \Delta''_{\omega,\,he^{-f}} & = & \Delta''_{\omega,\,h} + \bigg[[\Lambda_\omega,\,i\partial f\wedge\cdot\,],\,\bar\partial\bigg]:C^\infty_{p,\,q}(X,\,L')\longrightarrow {\cal D}'_{p,\,q}(X,\,L').\end{eqnarray}

\end{Def}

    In particular, the meaning of the identity on the second line in (\ref{eqn:Laplacians_singular}) is: \begin{eqnarray}\label{eqn:D'-star_singular}\bigg(D'_{he^{-f}}\bigg)^\star_{\omega,\,he^{-f}}\bigg(v^{(\alpha)}\otimes e_\alpha\bigg) = \bigg(\partial^\star_\omega v^{(\alpha)}\bigg)\otimes e_\alpha,\end{eqnarray} where $v^{(\alpha)}\otimes e_\alpha$ is the expression of an $L'$-valued $(p,\,q)$-form in a local trivialisation $\theta_\alpha:L'_{|U_\alpha}\longrightarrow U_\alpha\times\C$ of $L'$, $e_\alpha$ being the associated local frame of $L'$ and $v^{(\alpha)}$ being a $\C$-valued $(p,\,q)$-form on $U_\alpha$.

\vspace{1ex}    

    We now note that the classical Bochner-Kodaira-Nakano (BKN) identity, familiar from the smooth fibre-metric case, still holds in the singular setting. It involves the curvature current $i\Theta_{he^{-f}}(L') = i\Theta_h(L') + i\partial\bar\partial f$ of $(L',\,he^{-f})$, where $i\Theta_h(L')$ is the $C^\infty$ curvature $(1,\,1)$-form of $(L',\,h)$.

    \begin{Prop}\label{Prop:BKM_singular} In the setting of Definition \ref{Def:Laplacians_singular}, suppose that the metric $\omega$ on $X$ is {\bf K\"ahler}. Then, in every bidegree $(p,\,q)$, the following BKN identity of linear operators from $C^\infty_{p,\,q}(X,\,L')$ to ${\cal D}'_{p,\,q}(X,\,L')$ holds: \begin{eqnarray}\label{eqn:BKM_singular}\Delta''_{\omega,\,he^{-f}} = \Delta'_{\omega,\,he^{-f}} + \bigg[\bigg(i\Theta_h(L') + i\partial\bar\partial f\bigg)\wedge\cdot\,,\,\Lambda_\omega\bigg].\end{eqnarray}

\end{Prop}      
    
\noindent {\it Proof.} From Definition \ref{Def:Laplacians_singular} we get: \begin{eqnarray*}\Delta''_{\omega,\,he^{-f}} - \Delta'_{\omega,\,he^{-f}} = \Delta''_{\omega,\,h} - \Delta'_{\omega,\,h} + \bigg[[\Lambda_\omega,\,i\partial f\wedge\cdot\,],\,\bar\partial\bigg] + \bigg[[\Lambda_\omega,\,\bar\partial],\,i\partial f\wedge\cdot\,\bigg].\end{eqnarray*}

Now, the standard BKN identity applied to the $C^\infty$ fibre metric $h$ gives: \begin{eqnarray*}\Delta''_{\omega,\,h} - \Delta'_{\omega,\,h} = [i\Theta_h(L)\wedge\cdot\,,\,\Lambda_\omega],\end{eqnarray*} while the Jacobi identity spells: \begin{eqnarray*}\bigg[[\Lambda_\omega,\,i\partial f\wedge\cdot\,],\,\bar\partial\bigg] + \bigg[[i\partial f\wedge\cdot\,,\,\bar\partial],\,\Lambda_\omega\bigg] - \bigg[[\bar\partial,\,\Lambda_\omega],\,i\partial f\wedge\cdot\,\bigg] = 0.\end{eqnarray*}

Identity (\ref{eqn:BKM_singular}) follows from these equalities combined after noting that $[i\partial f\wedge\cdot\,,\,\bar\partial] = -i\partial\bar\partial f\wedge\cdot$ and $[\bar\partial,\,\Lambda_\omega] = - [\Lambda_\omega,\,\bar\partial]$. \hfill $\Box$

\subsubsection{Regularising $f$}\label{subsubsection:regularising} In the setting of Definition \ref{Def:Laplacians_singular}, let $m\in\{1,\dots , n\}$ and suppose that the curvature current $i\Theta_{he^{-f}}(L') = i\Theta_h(L') + i\partial\bar\partial f$ of $(L',\,he^{-f})$ is {\it $m$-semi-positive with respect to $\omega$}\footnote{This notion of partial positivity was introduced by Dieu in [Die06] and subsequently studied by several authors, including Harvey and Lawson in [HL13], Verbitsky in [Ver10], Dinew in [Din22], Dinew and Popovici in [DP25a] and [DP25b].}, a fact denoted by \begin{eqnarray}\label{eqn:m-pos_hypothesis}i\Theta_h(L') + i\partial\bar\partial f\geq_{m,\,\omega}0.\end{eqnarray} This means that the bidegree-$(m,\,m)$-current $i\Theta_{he^{-f}}(L')\wedge\omega^{m-1}$ is strongly semi-positive on $X$. Since the restriction of $i\Theta_h(L')$ to each open subset $U_\alpha\subset X$ on which $L$ is trivial equals $i\partial\bar\partial\varphi_\alpha$, hypothesis (\ref{eqn:m-pos_hypothesis}) means that, for every $\alpha$, the function $\varphi_\alpha + f$ is $m$-psh on $U_\alpha$ with respect to $\omega$.

By the local regularisation Theorem 1.6. in [DP25a], the open cover $(U_\alpha)_\alpha$ of $X$ can be chosen such that, for every $\alpha$, there exists a sequence $(f_{\alpha,\,\nu})_{\nu\geq 1}$ of $C^\infty$ functions $f_{\alpha,\,\nu}:U_\alpha\longrightarrow\R$ satisfying: \begin{eqnarray}\label{eqn:local-regularisations_prop}f_{\alpha,\,\nu}\searrow f_{|U_\alpha} \hspace{1ex} \mbox{pointwise on}\hspace{1ex} U_\alpha \hspace{1ex} \mbox{as}\hspace{1ex} \nu\to\infty \hspace{1ex}\mbox{and}\hspace{1ex} i\Theta_h(L')_{|U_\alpha} + i\partial\bar\partial f_{\alpha,\,\nu}>_{m,\,\omega}0 \hspace{1ex} \mbox{on}\hspace{1ex} U_\alpha \hspace{1ex} \mbox{for}\hspace{1ex} \nu\geq 1.\end{eqnarray}

Suppose that $X$ is compact and let $(\theta_\alpha)_\alpha$ be a smooth partition of unity on $X$ subordinate to the open cover $(U_\alpha)_\alpha$. For every $\nu\in\N^\star$, we define the following global $C^\infty$ function: \begin{eqnarray}\label{eqn:f_nu_def}f_\nu:=\sum\limits_\alpha\theta_\alpha f_{\alpha,\,\nu}:X\longrightarrow\R.\end{eqnarray} From (\ref{eqn:local-regularisations_prop}), we infer that $f_\nu\searrow f$ pointwise on $X$ as $\nu\to\infty$.


\subsection{Curvature of the $p$-contact and canonical bundles}\label{subsection:curvature_p-contact-canonical}

The proof of the main theorem in [Dem02] yields the following statement:

\begin{The}(Demailly)\label{The:Demailly_contact} Let $X$ be a compact {\bf K\"ahler} manifold with $\mbox{dim}_\C X = n$. Suppose there exist:

 \vspace{1ex} 

 $\bullet$ a Hermitian holomorphic line bundle $(L,\,h)$ on $X$, with possibly singular fibre metric $h=(e^{-\varphi_\alpha})_\alpha$, such that its dual $(L^{-1},\,h^{-1} = (e^{\varphi_\alpha})_\alpha)$ is {\bf pseudo-effective} in the sense that its curvature current $i\Theta_{h^{-1}}(L^{-1}) = -i\Theta_h(L)$ is $\geq 0$;

 \vspace{1ex} 

 $\bullet$ an integer $p\in\{0,\, 1,\dots , n\}$ and a non-zero $L$-valued holomorphic $(p,\,0)$-form $\Gamma = (\Gamma_\alpha)_\alpha\in H^{p,\,0}_{\bar\partial}(X,\,L)\setminus\{0\}$.

 \vspace{1ex}

 Then, $\partial\Gamma_\alpha = \partial\varphi_\alpha\wedge\Gamma_\alpha$ on $U_\alpha$ for every $\alpha$.

\end{The}

Since $i\Theta_h(L)_{|U_\alpha} = i\partial\bar\partial\varphi_\alpha$, the pseudo-effectivity assumption on $(L^{-1},\,h^{-1})$ means that each function $-\varphi_\alpha$ is plurisubharmonic on $U_\alpha$. In the special case where $n=2p+1$ and $p$ is odd, the conclusion means that $\Gamma = (\Gamma_\alpha)_\alpha\in H^{p,\,0}_{\bar\partial}(X,\,L)\setminus\{0\}$ is an {\bf $L$-valued holomorphic $p$-no-contact structure} on $X$ (see (\ref{eqn:L_hol-p-no-contact_Gamma_bis})). In particular, $\Gamma$ cannot be an $L$-valued holomorphic $p$-contact structure. We get the following consequence that was already expressed in a slightly different language in [Dem02].

\begin{Cor}\label{Cor:Demailly_K_not-psef} Let $X$ be a compact {\bf K\"ahler} manifold with $\mbox{dim}_\C X = n = 2p+1$ and $p$ odd.

\vspace{1ex}  

(i)\, If there exists a holomorphic line bundle $L$ on $X$ and an {\bf $L$-valued holomorphic $p$-contact structure} $\Gamma\in H^{p,\,0}_{\bar\partial}(X,\,L)$, then the dual line bundle $L^{-1}$, hence also the canonical line bundle $K_X$ of $X$, {\bf is not pseudo-effective}.

  In particular, the Kodaira dimension of $X$ is negative: $k(X) = -\infty$.

\vspace{1ex}  

\vspace{1ex}  

(ii)\, If, furthermore, the second Betti number of $X$ is $b_2(X) = 1$, then the anti-canonical bundle $-K_X$ is {\bf ample}, hence $X$ is a {\bf Fano} manifold.

\end{Cor}

The K\"ahler assumption is crucial in Demailly's theorem (whose proof uses integration by parts and would involve extra terms if the metric were not K\"ahler) and in the above statement. Recall that the $L$-valued holomorphic $p$-contact hypothesis on $\Gamma$ implies (see (i) of Proposition \ref{Prop:initial_L_p-contact}) that $L^{-1}$ is a square root of the canonical bundle $K_X$: \begin{eqnarray}\label{eqn:L-dual_square-root_K}(L^{-1})^2\simeq K_X.\end{eqnarray} Therefore, $L^{-1}$ is pseudo-effective if and only if $K_X$ is.

\vspace{1ex}

\noindent {\it End of proof of Corollary \ref{Cor:Demailly_K_not-psef}.} In view of Theorem \ref{The:Demailly_contact} and what was said above, only part (ii) still needs a proof. Since $X$ is compact K\"ahler, the Hodge decomposition and symmetry hold in every degree, in particular in degree $2$. Therefore, the extra hypothesis $b_2(X) = 1$ forces the equalities $h^{1,\,1}(X) = 1$ and $h^{2,\,0}(X) = h^{0,\,2}(X) = 0$. From this, we infer that every De Rham cohomology class of degree 2 is of pure type $(1,\,1)$, hence $H^2(X,\,\Z) = H^2(X,\,\Z)\cap H^{1,\,1}(X,\,\R)$ and this is of rank $1$. This implies that $X$ is projective and every holomorphic line bundle on $X$ is either positive (equivalently, ample), or flat, or negative. In particular, so is $K_X$.

However, $K_X$ cannot be either positive or flat since it would then be pseudo-effective, which is ruled out by part (i). It follows that $K_X$ must be negative (w.r.t. a $C^\infty$ fibre metric), which amounts to its dual $-K_X$ being positive (i.e. ample). This is equivalent to saying that $X$ is Fano. \hfill $\Box$

\vspace{3ex}

We will now strengthen a part of the conclusion of Corollary \ref{Cor:Demailly_K_not-psef} by means of the notion of $m$-semi-positivity introduced in [Die06] and the local regularisation result of [DP25a] (cf. $\S$\ref{subsubsection:regularising}).

\begin{The}\label{The:p-contact_n-p-pos} Let $L$ be a holomorphic line bundle on a compact {\bf K\"ahler} manifold $(X,\,\omega)$ with $\mbox{dim}_\C X = n = 2p+1$ and $p$ odd.

  If there exists an {\bf $L$-valued holomorphic $p$-contact structure} $\Gamma\in H^{p,\,0}_{\bar\partial}(X,\,L)$, there exists {\bf no} (singular or otherwise) Hermitian fibre metric $h^{-1}$ on the dual line bundle $L^{-1}$ (hence no such metric $h^{-2}$ on the canonical bundle $K_X\simeq L^{-2}$ of $X$) such that \begin{eqnarray*}i\Theta_{h^{-1}}(L^{-1})\geq_{n-p,\,\omega}0  \hspace{5ex} (\mbox{or equivalently} \hspace{3ex} i\Theta_{h^{-2}}(K_X)\geq_{n-p,\,\omega}0).\end{eqnarray*}

\end{The}

\noindent {\it Proof.} $\bullet$ To begin with, suppose that a $C^\infty$ Hermitian fibre metric $h$ exists on $L$ such that the induced fibre metric $h^{-1}$ on the dual line bundle $L^{-1}$ satisfies the partial semi-positivity property $i\Theta_{h^{-1}}(L^{-1})\geq_{n-p,\,\omega}0$. We will show that this leads to a contradiction.

Thanks to $\omega$ being K\"ahler, the well-known Bochner-Kodaira-Nakano (BKN) identity (cf. e.g. [Dem84]) applied to $(L,\,h)\longrightarrow (X,\,\omega)$ reads: \begin{eqnarray}\label{eqn:BKN_Kaehler}\Delta''_{\omega,\,h} = \Delta'_{\omega,\,h} + [i\Theta_h(L)\wedge\cdot\,,\,\Lambda_\omega].\end{eqnarray} This is an equality of operators acting on the $L$-valued $C^\infty$ forms of any bidegree on $X$.

Now, let $\lambda_1,\dots , \lambda_n:X\longrightarrow\R$ be the eigenvalues of the $(1,\,1)$-curvature form $i\Theta_h(L)$ of $(L,\,h)$ with respect to $\omega$. Meanwhile, let $x\in X$ be an arbitrary point and let $z_1,\dots , z_n$ be local holomorphic coordinates on $X$ centred at $x$ such that \begin{eqnarray}\label{eqn:coordinates_diagonalisation}\omega(x) = \sum\limits_{j=1}^n idz_j\wedge d\bar{z}_j  \hspace{5ex}\mbox{and}\hspace{5ex}  i\Theta_h(L)(x) = \sum\limits_{j=1}^n \lambda_j(x)\,idz_j\wedge d\bar{z}_j.\end{eqnarray}

It is standard (see e.g. [Dem97, VI. $\S5.2.$]) that, for every $L$-valued form $u=\sum\limits_{J,\,K}u_{J\overline{K}}\,dz_J\wedge d\bar{z}_K\otimes e_\alpha$ on an open subset $U_\alpha\subset X$ on which $L$ is trivialised by a non-vanishing section $e_\alpha$, one has: \begin{eqnarray}\label{eqn:curvature-operator_eigenvalues}[i\Theta_h(L)\wedge\cdot\,,\,\Lambda_\omega]\,u = \sum\limits_{J,\,K}\bigg(\sum\limits_{j\in J}\lambda_j + \sum\limits_{k\in K}\lambda_k - \sum\limits_{l=1}^n\lambda_l\bigg)\,u_{J\overline{K}}\,dz_J\wedge d\bar{z}_K\otimes e_\alpha    \hspace{7ex}\mbox{at}\hspace{1ex} x.\end{eqnarray}

In our case, when we apply this formula to the restriction $\Gamma_{|U_\alpha} = \sum\limits_{|J| = p}\Gamma_J^{(\alpha)}\,dz_J\otimes e_\alpha$ to $U_\alpha$ of the given holomorphic $p$-contact form $\Gamma$ in place of $u$, we get: \begin{eqnarray*}[i\Theta_h(L)\wedge\cdot\,,\,\Lambda_\omega]\,(\Gamma_{|U_\alpha}) = \sum\limits_{|J|=p}\bigg(\sum\limits_{j\in J}\lambda_j - \sum\limits_{l=1}^n\lambda_l\bigg)\,\Gamma^{(\alpha)}_J\,dz_J\otimes e_\alpha = - \sum\limits_{|J|=p}\lambda_{C_J}\Gamma_J^{(\alpha)}\,dz_J\otimes e_\alpha\end{eqnarray*} at $x$, where $C_J$ is the multi-index that is complementary to $J$ in $\{1,\dots , n\}$ and $\lambda_{C_J}$ is the sum of the eigenvalues $\lambda_j$ with $j\in C_J$. Hence, taking the pointwise inner product with $\Gamma_{|U_\alpha}$ and choosing $e_\alpha$ such that $|e_\alpha(x)|_h=1$, we get:  \begin{eqnarray*}\bigg\langle[i\Theta_h(L)\wedge\cdot\,,\,\Lambda_\omega]\,(\Gamma_{|U_\alpha}),\,\Gamma_{|U_\alpha}\bigg\rangle_{\omega,\,h} = -\sum\limits_{|J|=p}\lambda_{C_J}\,|\Gamma_J^{(\alpha)}|^2\geq 0  \hspace{7ex}\mbox{at}\hspace{1ex} x,\end{eqnarray*} where the last inequality holds if $-\lambda_{C_J}(x)\geq 0$ for every multi-index $J$ of length $|J|=p$.

  Since $i\Theta_{h^{-1}}(L^{-1}) = -i\Theta_h(L)$, the condition $-\lambda_{C_J}(x)\geq 0$ for every $J$ with $|J|=p$ is equivalent to the sum of any $n-p$ eigenvalues of the curvature form $i\Theta_{h^{-1}}(L^{-1})$ with respect to $\omega$ being non-negative at $x$. This is known to be equivalent to the condition $i\Theta_{h^{-1}}(L^{-1})\geq_{n-p,\,\omega}0$ at $x$.

  Thus, if the condition $i\Theta_{h^{-1}}(L^{-1})\geq_{n-p,\,\omega}0$ were satisfied at every point of $X$, we would get, after integrating the above pointwise inequality with respect to the volume form $dV_\omega:=\omega^n/n!$ induced by $\omega$, the inequality: \begin{eqnarray*}\bigg\langle\!\!\!\bigg\langle[i\Theta_h(L)\wedge\cdot\,,\,\Lambda_\omega]\,\Gamma,\, \Gamma\bigg\rangle\!\!\!\bigg\rangle_{\omega,\,h} \geq 0.\end{eqnarray*}

  When combined with the BKN identity (\ref{eqn:BKN_Kaehler}), with $\Delta''_{\omega,\,h}\Gamma = 0$ (due to $\bar\partial\Gamma = 0$, that we have by hypothesis, and to $\bar\partial^\star_{\omega,\,h}\Gamma = 0$, that we have for bidegree reasons) and with $\langle\langle\Delta'_{\omega,\,h}\Gamma,\,\Gamma\rangle\rangle_{\omega,\,h} = ||D'_h\Gamma||^2_{\omega,\,h} + ||(D'_h)^\star\Gamma||^2_{\omega,\,h} \geq 0$, the above inequality would imply $\langle\langle\Delta'_{\omega,\,h}\Gamma,\,\Gamma\rangle\rangle_{\omega,\,h} = 0$, hence also $D'_h\Gamma = 0$.

  Thus, the form $\Gamma\in H^{p,\,0}_{\bar\partial}(X,\,L)$ would be an $(L,\,h)$-valued holomorphic $p$-no-contact structure on $X$ (see Definition \ref{Def:local_hol-p-no-contact}), a fact that is incompatible with our holomorphic $p$-contact hypothesis on $\Gamma$ (as shown in a remark just after Definition \ref{Def:local_hol-p-no-contact}). This provides the sought-after contradiction and finishes the proof in the case of a $C^\infty$ Hermitian fibre metric.

\vspace{2ex}

$\bullet$ Suppose now that a singular Hermitian fibre metric $h$ exists on $L$ such that the induced fibre metric $h^{-1}$ on the dual line bundle $L^{-1}$ satisfies the partial semi-positivity property $i\Theta_{h^{-1}}(L^{-1})\geq_{n-p,\,\omega}0$. We will show that this too leads to a contradiction.

We will apply the construction described in $\S$\ref{subsection:Laplacians_singular}, in particular the local regularisation result of [DP25a, Theorem 1.6]. In our current notation, we let $h^{-1} = h_0\,e^{-f}$ for some $C^\infty$ Hermitian fibre metric $h_0$ on $L^{-1}$ and some upper semicontinuous $L^1_{loc}$ function $f:X\longrightarrow\R\cup\{-\infty\}$. Thus, our hypothesis reads: \begin{eqnarray*}i\Theta_{h^{-1}}(L^{-1}) = i\Theta_{h_0}(L^{-1}) + i\partial\bar\partial f\geq_{n-p,\,\omega}0.\end{eqnarray*} (Cf. (\ref{eqn:m-pos_hypothesis}) with $L'=L^{-1}$, $m=n-p$ and $h$ replaced by $h_0$.)

The BKN identity (\ref{eqn:BKM_singular}) for the line bundle $L$ equipped with the singular fibre metric $h_0^{-1}\,e^f$ reads: \begin{eqnarray*}\Delta''_{\omega,\,h_0^{-1}\,e^f} = \Delta'_{\omega,\,h_0^{-1}\,e^f} + \bigg[-\bigg(i\Theta_{h_0}(L^{-1}) + i\partial\bar\partial f\bigg)\wedge\cdot\,,\,\Lambda_\omega\bigg].\end{eqnarray*} Applying this to the holomorphic $p$-contact form $\Gamma\in H^{p,\,0}_{\bar\partial}(X,\,L)$, we get the second equality below: \begin{eqnarray*}\label{eqn:BKM_singular_proof} 0 = \int\limits_X\bigg\langle\Delta''_{\omega,\,h_0^{-1}\,e^f}\Gamma,\,\Gamma\bigg\rangle_{\omega,\,h_0^{-1}}\,e^f\,dV_\omega & = & \int\limits_X\bigg\langle\Delta'_{\omega,\,h_0^{-1}\,e^f}\Gamma,\,\Gamma\bigg\rangle_{\omega,\,h_0^{-1}}\,e^f\,dV_\omega \\
  \nonumber  & + & \int\limits_X\bigg\langle\bigg[-\bigg(i\Theta_{h_0}(L^{-1}) + i\partial\bar\partial f\bigg)\wedge\cdot\,,\,\Lambda_\omega\bigg]\,\Gamma,\,\Gamma\bigg\rangle_{\omega,\,h_0^{-1}}\,e^f\,dV_\omega,\end{eqnarray*} where the first equality follows from $\bar\partial\Gamma = 0$ and $\bar\partial^\star_{\omega,\,h_0^{-1}\,e^f}\Gamma = 0$ (or, equivalently, from the definition of $\Delta''_{\omega,\,h_0^{-1}\,e^f}$ for a singular $f$ given on the last line of (\ref{eqn:Laplacians_singular}) and the fact that the $(p,\,0)$-form $\Gamma$ lies in the kernel of $\bar\partial$ and of any operator of bidegree $(0,\,-1)$). Note that the well-definedness of the integrals in the second equality above (for example, that of the current $e^f\,\partial\bar\partial f$ that features in the last integral) follows from (\ref{eqn:weak-convergences}) and the discussion immediately thereafter.   

We will show that \begin{eqnarray}\label{eqn:Delta'_singular_parts}\nonumber\int\limits_X\bigg\langle\Delta'_{\omega,\,h_0^{-1}e^f}\Gamma,\,\Gamma\bigg\rangle_{\omega,\,h_0^{-1}}e^f\,dV_\omega & = & \int\limits_X\bigg|D'_{h_0^{-1}e^f}\Gamma\bigg|^2_{\omega,\,h_0^{-1}}\,e^f\,dV_\omega + \int\limits_X\bigg|\bigg(D'_{h_0^{-1}e^f}\bigg)^\star_{\omega,\,h_0^{-1}e^f}\Gamma\bigg|^2_{\omega,\,h_0^{-1}}\,e^f\,dV_\omega \\
 & \geq & 0 \end{eqnarray} (where the meaning of the adjoint operator in the last integrand is spelt out in (\ref{eqn:D'-star_singular})) and \begin{eqnarray}\label{eqn:curvature-op_singular_pos}\int\limits_X\bigg\langle\bigg[-\bigg(i\Theta_{h_0}(L^{-1}) + i\partial\bar\partial f\bigg)\wedge\cdot\,,\,\Lambda_\omega\bigg]\,\Gamma,\,\Gamma\bigg\rangle_{\omega,\,h_0^{-1}}\,e^f\,dV_\omega\geq 0,\end{eqnarray} which will imply that $D'_{h_0^{-1}e^f}\Gamma = 0$. As in the case of a $C^\infty$ fibre metric, this will mean that $\partial\Gamma_\alpha = \partial(\psi_\alpha - f)\wedge\Gamma_\alpha$ on $U_\alpha$ for every $\alpha$ (where we put $h_0^{-1} = (e^{-\psi_\alpha})_\alpha$), which will lead to $\Gamma_\alpha\wedge\partial\Gamma_\alpha = -\partial(\psi_\alpha - f)\wedge\Gamma_\alpha\wedge\Gamma_\alpha = 0$, contradicting the hypothesis $\Gamma_\alpha\wedge\partial\Gamma_\alpha\neq 0$ everywhere on $U_\alpha$.

\vspace{1ex}

The arguments of $\S$\ref{subsubsection:regularising} yield a sequence $(f_\nu)_\nu$ of $C^\infty$ functions $f_\nu:X\longrightarrow\R$ such that $f_\nu\searrow f$ pointwise on $X$ as $\nu\to\infty$.

In particular, for every constant $c>0$, the functions $(e^{cf_\nu})_{\nu\geq 1}$ are uniformly bounded. Through a standard argument (using that $e^{cf_\nu}\leq e^{cf_1}\in L^1_{loc}$ for all $\nu\geq 1$, the dominated convergence theorem and the distributional derivatives $\partial(e^{cf})$, $\bar\partial(e^{cf})$, $\partial\bar\partial(e^{cf})$), we infer the following convergences in the weak topology of currents as $\nu\to\infty$: \begin{eqnarray}\label{eqn:weak-convergences}e^{cf_\nu}\longrightarrow e^{cf}, \hspace{3ex} e^{cf_\nu}\partial f_\nu\longrightarrow e^{cf}\partial f, \hspace{3ex} e^{cf_\nu}\partial f_\nu\wedge\bar\partial f_\nu\longrightarrow e^{cf}\partial f\wedge\bar\partial f, \hspace{3ex} e^{cf_\nu}\partial\bar\partial f_\nu\longrightarrow e^{cf}\partial\bar\partial f.\end{eqnarray} Note that the current $e^{cf}\partial f$ is well defined as the distributional $\partial$-derivative $\frac{1}{c}\,\partial(e^{c f})$; the current $e^{cf}\partial f\wedge\bar\partial f$ is defined as the weak limit of $e^{cf_\nu}\partial f_\nu\wedge\bar\partial f_\nu$ as $\nu\to\infty$; finally, the current $e^{cf}\partial\bar\partial f$ is defined as the difference $(1/c)\,\partial\bar\partial(e^{cf}) - c\,e^{cf}\partial f\wedge\bar\partial f$ of bona fide currents. (Alternatively, one can first define the current $e^{cf}\partial\bar\partial f$ as the weak limit of $e^{cf_\nu}\partial\bar\partial f_\nu$ as $\nu\to\infty$ and then the current $e^{cf}\partial f\wedge\bar\partial f$ as the difference $(1/c)\,\partial\bar\partial(e^{cf}) - e^{cf}\partial\bar\partial f$ of bona fide currents. The result is the same as above and the currents defined as weak limits are independent of the choices of smooth approximations.) Further note that this procedure can be carried out since in the equalities: \begin{eqnarray*}\frac{1}{c}\,i\partial\bar\partial(e^{cf_{\alpha,\,\nu}}) = c\,e^{cf_{\alpha,\,\nu}}\,i\partial f_{\alpha,\,\nu}\wedge\bar\partial f_{\alpha,\,\nu} + e^{cf_{\alpha,\,\nu}}\,i\partial\bar\partial f_{\alpha,\,\nu}  \hspace{3ex} \mbox{on}\hspace{1ex} U_\alpha, \hspace{3ex} \nu\geq 1,\end{eqnarray*} the left-hand-side term is bounded in mass independently of $\nu\geq 1$ (because it converges weakly, as $\nu\to\infty$, to the distributional $\partial\bar\partial$-derivative $(1/c)\,i\partial\bar\partial(e^{cf})$) and both of the right-hand-side terms have non-negative masses (because $e^{cf_{\alpha,\,\nu}}\,i\partial f_{\alpha,\,\nu}\wedge\bar\partial f_{\alpha,\,\nu}\wedge\omega^{n-1}\geq 0$ and $e^{cf_{\alpha,\,\nu}}\,i\partial\bar\partial f_{\alpha,\,\nu}\wedge\omega^{n-1}\geq 0)$, so their masses are bounded above by the one of the left-hand-side term and, therefore, they, too, are bounded independently of $\nu\geq 1$.

\vspace{1ex}  

For every $\nu\geq 1$, the standard BKN identity applied to the $C^\infty$ fibre metric $h_0^{-1}\,e^{f_\nu}$ on $L$ yields, after integrating and expanding the adjoints, the analogue of (\ref{eqn:Delta'_singular_parts}) with $f_\nu$ in place of $f$: \begin{eqnarray*}\bigg\langle\bigg\langle\Delta'_{\omega,\,h_0^{-1}\,e^{f_\nu}}\Gamma,\,\Gamma\bigg\rangle\bigg\rangle_{\omega,\,h_0^{-1}\,e^{f_\nu}} = \bigg|\bigg|D'_{h_0^{-1}\,e^{f_\nu}}\Gamma\bigg|\bigg|^2_{\omega,\,h_0^{-1}\,e^{f_\nu}} + \bigg|\bigg|\bigg(D'_{h_0^{-1}\,e^{f_\nu}}\bigg)^\star_{\omega,\,h_0^{-1}\,e^{f_\nu}}\Gamma\bigg|\bigg|^2_{\omega,\,h_0^{-1}\,e^{f_\nu}}.\end{eqnarray*} Given the expressions (\ref{eqn:Laplacians_singular}) and (\ref{eqn:D'-star_singular}), which hold for the operators in question when the fibre metric is either $C^\infty$ or singular, this translates to: \begin{eqnarray*} & & \int\limits_X\bigg\langle\Delta'_{\omega,\,h_0^{-1}}\Gamma,\,\Gamma\bigg\rangle_{\omega,\,h_0^{-1}}e^{f_\nu}\,dV_\omega + \int\limits_X\bigg\langle\bigg[[\Lambda_\omega,\,\bar\partial],\,i\partial f_\nu\wedge\cdot\,\bigg]\,\Gamma,\,\Gamma\bigg\rangle_{\omega,\,h_0^{-1}}\,e^{f_\nu}\,dV_\omega  \\
  &  = & \int\limits_X\bigg|\partial\Gamma-\partial\psi_\alpha\wedge\Gamma\bigg|^2_{\omega,\,h_0^{-1}}e^{f_\nu}\,dV_\omega + \int\limits_X\bigg|e^{f_\nu/2}\,\partial f_\nu\wedge\Gamma\bigg|^2_{\omega,\,h_0^{-1}}\,dV_\omega \\
 & + & \int\limits_X\bigg\langle\partial\Gamma-\partial\psi_\alpha\wedge\Gamma,\,e^{f_\nu}\partial f_\nu\wedge\Gamma\bigg\rangle_{\omega,\,h_0^{-1}}\,dV_\omega + \int\limits_X\bigg\langle e^{f_\nu}\partial f_\nu\wedge\Gamma,\,\partial\Gamma-\partial\psi_\alpha\wedge\Gamma\bigg\rangle_{\omega,\,h_0^{-1}}\,dV_\omega  \\
    &  + & \int\limits_X\bigg|\partial^\star_\omega\Gamma\bigg|^2_{\omega,\,h_0^{-1}}e^{f_\nu}\,dV_\omega ,\end{eqnarray*} where we set $h_0^{-1}:=(e^{-\psi_\alpha})_\alpha$. The weak convergences (\ref{eqn:weak-convergences}) imply that each term in the above equality converges, as $\nu\to\infty$, to its analogue with $f$ in place of $f_\nu$. Thanks to the expressions (\ref{eqn:Laplacians_singular}) for the operators in question induced by the singular fibre metric $h_0^{-1}\,e^f$ on $L$, this means that (\ref{eqn:Delta'_singular_parts}) holds.

  \vspace{1ex}

  To prove inequality (\ref{eqn:curvature-op_singular_pos}), recall from $\S$\ref{subsubsection:regularising} that $f_\nu:=\sum\limits_\alpha\theta_\alpha f_{\alpha,\,\nu}$ with $C^\infty$ functions $f_{\alpha,\,\nu}:U_\alpha\to\R$ satisfying $f_{\alpha,\,\nu}\searrow f_{|U_\alpha}$ pointwise on $U_\alpha$ as $\nu\to\infty$ and \begin{eqnarray*}i\Theta_{h_0}(L^{-1})_{|U_\alpha} + i\partial\bar\partial f_{\alpha,\,\nu}>_{n-p,\,\omega}0 \hspace{1ex} \mbox{on}\hspace{1ex} U_\alpha \hspace{1ex} \mbox{for}\hspace{1ex} \nu\geq 1.\end{eqnarray*} As shown in the first part of the proof (for a $C^\infty$ fibre metric), this positivity property implies: \begin{eqnarray*}\bigg\langle\bigg[-\bigg(i\Theta_{h_0}(L^{-1})_{|U_\alpha} + i\partial\bar\partial f_{\alpha,\,\nu}\bigg)\wedge\cdot\,,\,\Lambda_\omega\bigg]\,(\Gamma_{|U_\alpha}),\,(\Gamma_{|U_\alpha})\bigg\rangle_{\omega,\,h_0^{-1}}\,e^{f_{\alpha,\,\nu}}\geq 0\end{eqnarray*} at every point of $U_\alpha$ for every $\nu\geq 1$. Meanwhile, thanks to the weak convergences (\ref{eqn:weak-convergences}), the left-hand side of the above inequality converges weakly, as $\nu\to\infty$, to the integrand on the left-hand side of (\ref{eqn:curvature-op_singular_pos}). This proves inequality (\ref{eqn:curvature-op_singular_pos}) and we are done.  \hfill $\Box$

\subsection{Curvature of the contact and canonical bundles}\label{subsection:curvature_contact-canonical}

By essentially the same argument, we can show that the existence of a holomorphic contact structure rules out the canonical bundle $K_X$ being partially semi-positive in an even weaker sense.

\begin{The}\label{The:contact_n-1-pos} Let $L$ be a holomorphic line bundle on a compact {\bf K\"ahler} manifold $(X,\,\omega)$ with $\mbox{dim}_\C X = n = 2p+1$. 

  If there exists an {\bf $L$-valued holomorphic contact structure} $\eta\in H^{1,\,0}_{\bar\partial}(X,\,L)$, there exists {\bf no} (singular or otherwise) Hermitian fibre metric $h^{-1}$ on the dual line bundle $L^{-1}$ (hence no such metric $h^{-(p+1)}$ on the canonical bundle $K_X\simeq L^{-(p+1)}$ of $X$) such that \begin{eqnarray*}i\Theta_{h^{-1}}(L^{-1})\geq_{n-1,\,\omega}0  \hspace{5ex} (\mbox{or equivalently} \hspace{3ex} i\Theta_{h^{-(p+1)}}(K_X)\geq_{n-1,\,\omega}0).\end{eqnarray*}

\end{The}

\noindent {\it Proof.} We will only point out the differences to the proof of Theorem \ref{The:p-contact_n-p-pos}. We keep the same notation. When applied to the restriction $\eta_{|U_\alpha} = \sum\limits_{j = 1}^n\eta_j^{(\alpha)}\,dz_j\otimes e_\alpha$, the general formula (\ref{eqn:curvature-operator_eigenvalues}) yields: \begin{eqnarray*}[i\Theta_h(L)\wedge\cdot\,,\,\Lambda_\omega]\,(\eta_{|U_\alpha}) = \sum\limits_{j=1}^n\bigg(\lambda_j - \sum\limits_{l=1}^n\lambda_l\bigg)\,\eta^{(\alpha)}_j\,dz_j\otimes e_\alpha = - \sum\limits_{j=1}^n\lambda_{C_j}\eta_j^{(\alpha)}\,dz_j\otimes e_\alpha\end{eqnarray*} at $x$, where $\lambda_{C_j}:=\sum_{l=1}^n\lambda_l - \lambda_j$. If we choose $e_\alpha$ such that $|e_\alpha(x)|_h=1$, it follows that \begin{eqnarray*}\bigg\langle[i\Theta_h(L)\wedge\cdot\,,\,\Lambda_\omega]\,(\eta_{|U_\alpha}),\,\eta_{|U_\alpha}\bigg\rangle_{\omega,\,h} = -\sum\limits_{j=1}^n\lambda_{C_j}\,|\eta_j^{(\alpha)}|^2\geq 0  \hspace{7ex}\mbox{at}\hspace{1ex} x,\end{eqnarray*} where the last inequality holds if $-\lambda_{C_j}(x)\geq 0$ for every $j\in\{1,\dots , n\}$.

  Now, the condition $-\lambda_{C_j}(x)\geq 0$ for every $j\in\{1,\dots , n\}$ is equivalent to the sum of any $n-1$ eigenvalues with respect to $\omega$ of the curvature form $i\Theta_h(L)$ of $(L,\,h)$ being non-positive at $x$, hence to  the sum of any $n-1$ eigenvalues of the curvature form $i\Theta_{h^{-1}}(L^{-1})$ being non-negative at $x$. This is further equivalent to the condition $i\Theta_{h^{-1}}(L^{-1})\geq_{n-1,\,\omega}0$ at $x$.

  If $L$ has a $C^\infty$ fibre metric $h$ satisfying this property, we get a contradiction to $\eta$ being a holomorphic $L$-valued contact structure on $X$ in the same way as in the proof of Theorem \ref{The:p-contact_n-p-pos}. Specifically, this would imply that $D_h'\eta = 0$, a fact that amounts to $\partial\eta_\alpha = \partial\varphi_\alpha\wedge\eta_\alpha$ on $U_\alpha$ for every $\alpha$. Then, we would have $\eta_\alpha\wedge\partial\eta_\alpha = 0$ (since $\eta_\alpha$ is a $1$-form, so $\eta_\alpha\wedge\eta_\alpha = 0$), hence also $\eta_\alpha\wedge(\partial\eta_\alpha)^p = 0$, in violation of the holomorphic contact hypothesis made on $\eta$ (cf. (\ref{eqn:L_hol-contact_eta})).

 If $L$ has a singular fibre metric $h$ satisfying the curvature hypothesis, we locally regularise it using Theorem 1.6 in [DP25a] and then we argue via the construction described in $\S$\ref{subsection:Laplacians_singular} as we did in the proof of Theorem \ref{The:p-contact_n-p-pos}. \hfill $\Box$

\subsection{Scalar curvature of spin manifolds}\label{subsection:scalar-spin}

We shall now see that, moving in the opposite direction towards the weaker notion of spin manifold, the same type of argument yields a vanishing theorem under a positivity assumption on the scalar curvature of the manifold. If $\omega$ is a Hermitian metric on an $n$-dimensional complex manifold $X$ and $\tilde{h}$ is a $C^\infty$ Hermitian fibre metric on the anti-canonical bundle $-K_X$, we denote by \begin{eqnarray*}{\mbox Scal}\,(\omega,\,\tilde{h}):=\Lambda_\omega\bigg(\frac{i}{2\pi}\,\Theta_{\tilde{h}}(-K_X)\bigg)\end{eqnarray*} the {\it scalar curvature} of $X$ with respect to $\omega$ and $\tilde{h}$. This is a real-valued function on $X$ that generalises the usual scalar curvature function ${\mbox Scal}\,(\omega):=\Lambda_\omega({\mbox Ric}\,\omega)$, the latter being obtained from ${\mbox Scal}\,(\omega,\,\tilde{h})$ by taking $\tilde{h}$ to be the $C^\infty$ Hermitian fibre metric on $-K_X$ induced by $\omega$.

\begin{The}\label{The:spin_scal-pos} Let $X$ be a compact connected complex manifold with $\mbox{dim}_\C X = n$. Suppose there exists a {\bf spin structure} on $X$, namely a holomorphic line bundle $L$ on $X$ such that $L^2 = -K_X$.

  If there exists a {\bf K\"ahler} metric $\omega$ on $X$ and a $C^\infty$ Hermitian fibre metric $h$ on $L$ (equivalently, a $C^\infty$ Hermitian fibre metric $h^2$ on $-K_X$) such that \begin{eqnarray}\label{eqn:Scal_positivity_point}{\mbox Scal}\,(\omega,\,h^2)\geq 0 \hspace{2ex}\mbox{everywhere on}\hspace{1ex} X \hspace{3ex}\mbox{and}\hspace{3ex} {\mbox Scal}\,(\omega,\,h^2)> 0 \hspace{2ex}\mbox{at some point of}\hspace{1ex} X,\end{eqnarray} then $H^0(X,\,L^{-1}) = \{0\}$ and $H^0(X,\,K_X) = H^{n,\,0}(X,\,\C) = \{0\}$.

\end{The}   

\noindent {\it Proof.} The notation is carried forward from the proofs of Theorems \ref{The:p-contact_n-p-pos} and \ref{The:contact_n-1-pos}. We fix a K\"ahler metric $\omega$ on $X$ and a $C^\infty$ Hermitian fibre metric $h$ on $L$.

The general Bochner-Kodaira-Nakano (BKN) identity (\ref{eqn:BKN_Kaehler}) applied to $(L^{-1},\,h^{-1})\longrightarrow (X,\,\omega)$ in any bidegree reads: \begin{eqnarray*}\Delta''_{\omega,\,h^{-1}} = \Delta'_{\omega,\,h^{-1}} + [i\Theta_{h^{-1}}(L^{-1})\wedge\cdot\,,\,\Lambda_\omega].\end{eqnarray*} Therefore, in bidegree $(0,\,0)$, for every global smooth section $f\in C^\infty_{0,\,0}(X,\,L^{-1})$, we get after taking $L^2_{\omega,\,h^{-1}}$-inner products on $X$: \begin{eqnarray}\label{eqn:Scal_positivity_point_proof_1}\nonumber\bigg\langle\!\!\!\bigg\langle\Delta''_{\omega,\,h^{-1}}f,\,f\bigg\rangle\!\!\!\bigg\rangle_{\omega,\,h^{-1}} & = & \bigg\langle\!\!\!\bigg\langle\Delta'_{\omega,\,h^{-1}}f,\,f\bigg\rangle\!\!\!\bigg\rangle_{\omega,\,h^{-1}} + \bigg\langle\!\!\!\bigg\langle[i\Theta_{h^{-1}}(L^{-1})\wedge\cdot\,,\,\Lambda_\omega]\,f,\,f\bigg\rangle\!\!\!\bigg\rangle_{\omega,\,h^{-1}} \\
  & \geq & \int\limits_X \bigg\langle[i\Theta_{h^{-1}}(L^{-1})\wedge\cdot\,,\,\Lambda_\omega]\,f,\,f\bigg\rangle_{h^{-1}}\,dV_\omega.\end{eqnarray}

Meanwhile, after writing $f_{|U_\alpha} = f_\alpha\otimes e_\alpha$ on $U_\alpha$, the general formula (\ref{eqn:curvature-operator_eigenvalues}) yields: \begin{eqnarray*}[i\Theta_{h^{-1}}(L^{-1})\wedge\cdot\,,\,\Lambda_\omega]\,(f_{|U_\alpha}) = -\sum\limits_{j=1}^n(-\lambda_j)\,f^{(\alpha)}\otimes e_\alpha   \hspace{7ex}\mbox{at}\hspace{1ex} x.\end{eqnarray*} Since $\sum\limits_{j=1}^n\lambda_j = \Lambda_\omega(i\Theta_h(L)) = (1/2)\,\Lambda_\omega(i\Theta_{h^2}(-K_X)) = \pi\,{\mbox Scal}\,(\omega,\,h^2)$ at $x$, we get: \begin{eqnarray}\label{eqn:Scal_positivity_point_proof_2}\bigg\langle[i\Theta_{h^{-1}}(L^{-1})\wedge\cdot\,,\,\Lambda_\omega]\,f,\,f\bigg\rangle_{h^{-1}} = \pi\,{\mbox Scal}\,(\omega,\,h^2)\,|f|^2_h\geq 0 \hspace{7ex}\mbox{everywhere on}\hspace{1ex} X\end{eqnarray} since $x\in X$ was arbitrary. The last inequality above follows from our hypothesis (\ref{eqn:Scal_positivity_point}).

Suppose now that $f$ is holomorphic (i.e. $f\in H^0(X,\,L^{-1})$). Then, $\Delta''_{\omega,\,h^{-1}}f = 0$, so the combined (\ref{eqn:Scal_positivity_point_proof_1}) and (\ref{eqn:Scal_positivity_point_proof_2}) imply that ${\mbox Scal}\,(\omega,\,h^2)\,|f|^2_h = 0$ at every point of $X$. Since the function ${\mbox Scal}\,(\omega,\,h^2)$ is continuous (even $C^\infty$) on $X$ and positive at some point (by hypothesis (\ref{eqn:Scal_positivity_point})), it must be positive on a neighbourhood of that point. But then $f$ must vanish identically on that neighbourhood. Since $f$ is holomorphic and $X$ is connected, we infer that $f$ vanishes at every point of $X$. \hfill $\Box$

\vspace{3ex}

One of the consequences of this discussion is the following

\begin{Cor}\label{Cor:p-contact_non-Kobayashi-hyperbolic} Let $X$ be a projective manifold with $\mbox{dim}_\C X = n = 2p+1$ and $p$ odd. If there exists a {\bf holomorphic $p$-contact structure} $\Gamma\in H^{p,\,0}_{\bar\partial}(X,\,L)$ with values in a holomorphic line bundle $L$ over $X$, then $X$ is {\bf covered by rational curves}. In particular, $X$ is {\bf not Kobayashi hyperbolic}.

\end{Cor}

\noindent {\it Proof.} By Corollary \ref{Cor:Demailly_K_not-psef}, the canonical line bundle $K_X$ is not pseudo-effective. Then, by Corollary 0.3. in [BDPP13] (a consequence of the main cone duality Theorem 0.2. in [BDPP13] and of the classical Miyaoka-Mori bend-and-break lemma in characteristic $p$), $X$ is covered by rational curves.

Meanwhile, the existence in $X$ of a rational curve precludes the Kobayashi hyperbolicity of $X$, as is well known and easy to check.  \hfill $\Box$

\section{Local holomorphic $s$-symplectic structures}\label{section:local_hol-s-symplecticcontact}

We now introduce our second main notion: the local (i.e. line bundle-valued) analogue of the notion of scalar holomorphic $s$-symplectic structure introduced in [KPU25]. 

\begin{Def}\label{Def:local_hol-s-symplectic} Let $L$ be a holomorphic line bundle on a compact complex manifold $X$ with $\mbox{dim}_\C X = n = 2s$, where $s$ is an even integer.

 \vspace{1ex}

 $(1)$\, An {\bf $L$-valued holomorphic $s$-symplectic structure} on $X$ is a form $\Omega\in C^\infty_{s,\,0}(X,\,L)$ satisfying the following two conditions: \begin{eqnarray}\label{eqn:L_hol-s-symplectic_Omega}(a)\,\,\, \bar\partial\Omega = 0 \hspace{5ex}\mbox{and}\hspace{5ex} (b)\,\,\,\Omega\wedge\Omega \neq 0 \hspace{2ex}\mbox{at every point of}\hspace{2ex} X.\end{eqnarray}

 When the line bundle $L$ is not specified, any such form $\Omega$ is alternatively called a {\bf local holomorphic $s$-symplectic structure} on $X$.

 \vspace{1ex}

$(2)$\, If there exists an $L$-valued holomorphic $s$-symplectic structure $\Omega\in C^\infty_{s,\,0}(X,\,L)$, $L$ is called a {\bf holomorphic $s$-symplectic bundle}, $(X,\,L,\,\Omega)$ is called a {\bf holomorphic $s$-symplectic triple} and $X$ is called a {\bf (locally holomorphic) $s$-symplectic manifold}.

\end{Def}

Note that this definition is correct since $\Omega\wedge\Omega$ is a well-defined $C^\infty$ (even holomorphic) $(n,\,0)$-form on $X$ with values in $L^2$. Indeed, if $\Omega = (\Omega_\alpha)_\alpha$ with respect to a family $(U_\alpha,\,\theta_\alpha,\,e_\alpha)_\alpha$ of local holomorphic trivialisations of $L$, we have $\Omega_\alpha = g_{\alpha\beta}\,\Omega_\beta$, which implies that \begin{eqnarray*}\Omega_\alpha\wedge\Omega_\alpha =  g_{\alpha\beta}^2\,\Omega_\beta\wedge\Omega_\beta \hspace{5ex} \mbox{on}\hspace{2ex} U_\alpha\cap U_\beta \hspace{2ex} \mbox{for all} \hspace{2ex}\alpha,\beta.\end{eqnarray*} Note that for condition (b) in (\ref{eqn:L_hol-s-symplectic_Omega}) to be satisfied, it is necessary that $s$ be {\it even} since otherwise $\Omega_\alpha\wedge\Omega_\alpha$ would be the zero form on $U_\alpha$ for every $\alpha$. 

Thus, if $\Omega = (\Omega_\alpha)_\alpha$ is an $L$-valued holomorphic $s$-symplectic structure on $X$, then \begin{eqnarray*}\Omega\wedge\Omega = (\Omega_\alpha\wedge\Omega_\alpha)_\alpha\in   H^{n,\,0}(X,\,L^2)\simeq H^0(X,\,\Omega^n_X\otimes L^2) = H^0(X,\,K_X\otimes L^2)\end{eqnarray*} is a {\bf non-vanishing} global holomorphic section of the holomorphic line bundle $K_X\otimes L^2$. This line bundle is, therefore, necessarily holomorphically trivial, hence we get a holomorphic line bundle isomorphism $L^2\simeq -K_X$. In particular, $-K_X$ and $K_X$ have square roots, so the classical Proposition \ref{Prop:spin_existence_standard} yields the following analogue in this context of Corollary \ref{Cor:p-contact_spin}.

\begin{Cor}\label{Cor:s-symplectic_spin} Any compact complex manifold $X$ with $\mbox{dim}_\C X \equiv 0 \hspace{2ex} mod \hspace{1ex} 4$ that carries a {\bf local holomorphic $s$-symplectic structure} is a {\bf spin manifold}.

\end{Cor}

Similarly to the case of holomorphic $p$-contact structures, we get the following

\begin{Conc}\label{Conc:all-p-contact_L} Let $L$ be a holomorphic line bundle over a compact complex manifold $X$ with $\mbox{dim}_\C X = 2s$ and $s$ even such that $L^2\simeq -K_X$.

  Then, the set of $L$-valued holomorphic $s$-symplectic structures $\Omega$ on $X$ is either empty or the complement of a degree-$2$ complex hypersurface in the complex projective space $\Proj(H^{s,\,0}_{\bar\partial}(X,\,L))$.

\end{Conc}

This follows by considering the entire holomorphic map \begin{eqnarray*}\label{eqn:Omega-map-homogeneous}H^{s,\,0}_{\bar\partial}(X,\,L)\ni\Omega\stackrel{S}{\longmapsto}\Omega\wedge\Omega\in\C,\end{eqnarray*} observing that it is homogeneous of degree $2$ and concluding that it must be a homogeneous polynomial of degree $2$, hence it identifies with a global holomorphic section of ${\cal O}(2)$ over $\Proj(H^{s,\,0}_{\bar\partial}(X,\,L))$. The $L$-valued holomorphic $s$-symplectic structures $\Omega$ on $X$, if any, are precisely the elements $\Omega\in H^{s,\,0}_{\bar\partial}(X,\,L)$ such that $S(\Omega)\neq 0$.

\vspace{3ex}

One link between $p$-contact and $s$-symplectic manifolds is provided by the following

\begin{Prop}\label{Prop:p-contact_s-symplectic_fibration} Let $\pi:X\longrightarrow Y$ be a surjective holomorphic submersion between compact complex manifolds with $\mbox{dim}_\C Y = 2s = 4r$ and $\mbox{dim}_\C X = 2s + 2p + 1$, where $s,r,p$ are positive integers and $p$ is odd. Suppose there exist:


  $\bullet$ a holomorphic line bundle $L$ on $Y$ and an $L$-valued holomorphic {\bf $s$-symplectic structure} $\Omega\in C^\infty_{s,\,0}(Y,\,L)$;


  $\bullet$ a holomorphic line bundle $F$ on $X$ and a form $\Gamma\in C^\infty_{p,\,0}(X,\,F)$ such that $\bar\partial\Gamma = 0$ on $X$ and, for every $y\in Y$, $\Gamma_y\wedge\partial\Gamma_y\neq 0$ everywhere on $X_y:=\pi^{-1}(y)$, where $\Gamma_y:=\Gamma_{|X_y}$.

  Then, the form $\widetilde\Gamma := \pi^\star\Omega\wedge\Gamma\in C^\infty_{s+p,\,0}(X,\,\pi^\star L\otimes F)$ is a $(\pi^\star L\otimes F)$-valued holomorphic {\bf $(s+p)$-contact structure} on $X$.

\end{Prop}  

\noindent {\it Proof.} Since $\pi$ is holomorphic, $\bar\partial\Omega = 0$ and  $\bar\partial\Gamma = 0$, we get: \begin{eqnarray*}\bar\partial\widetilde\Gamma = \pi^\star(\bar\partial\Omega)\wedge\Gamma + \pi^\star\Omega\wedge\bar\partial\Gamma = 0.\end{eqnarray*}

On the other hand, we get: \begin{eqnarray*}\widetilde\Gamma\wedge\partial\widetilde\Gamma = \pi^\star\Omega\wedge\Gamma\wedge\bigg(\partial(\pi^\star\Omega)\wedge\Gamma + \pi^\star\Omega\wedge\partial\Gamma\bigg) = \pi^\star(\Omega\wedge\Omega)\wedge(\Gamma\wedge\partial\Gamma),\end{eqnarray*} where the last identity follows from $\Gamma\wedge\Gamma = 0$, itself due to $p$ being odd. The hypotheses imply that this holomorphic top-degree form on $X$ is nowhere vanishing. \hfill $\Box$

\vspace{2ex}

The $s$-symplectic analogue of Theorem \ref{The:p-contact_n-p-pos} is the following

\begin{The}\label{The:s-symplectic_s-pos} Let $L$ be a holomorphic line bundle on a compact {\bf K\"ahler} manifold $(X,\,\omega)$ with $\mbox{dim}_\C X = n = 2s$ and $s$ even. Suppose there exist an {\bf $L$-valued holomorphic $s$-symplectic structure} $\Omega\in H^{s,\,0}_{\bar\partial}(X,\,L)$ and a (possibly singular) Hermitian fibre metric $h^{-1}$ on the dual line bundle $L^{-1}$ (or, equivalently, such a metric $h^{-2}$ on the canonical bundle $K_X\simeq L^{-2}$ of $X$) such that \begin{eqnarray*}i\Theta_{h^{-1}}(L^{-1})\geq_{s,\,\omega}0  \hspace{5ex} (\mbox{or, equivalently,} \hspace{3ex} i\Theta_{h^{-2}}(K_X)\geq_{s,\,\omega}0).\end{eqnarray*}

  Then, the sheaf ${\cal F}_\Omega$ of germs of holomorphic $(1,\,0)$-vector fields $\xi$ on $X$ such that $\xi\lrcorner\Omega = 0$ has rank zero.

\end{The}

\noindent {\it Proof.} $\bullet$ Suppose there exists an $L$-valued holomorphic $s$-symplectic structure $\Omega\in H^{s,\,0}_{\bar\partial}(X,\,L)$. If $h$ is a $C^\infty$ Hermitian fibre metric on $L$ and $\omega$ is a K\"ahler metric on $X$, the Bochner-Kodaira-Nakano identity (\ref{eqn:BKN_Kaehler}) yields the second equality below: \begin{eqnarray*}0 = \bigg\langle\!\!\!\bigg\langle\Delta''_{\omega,\,h}\Omega,\,\Omega\bigg\rangle\!\!\!\bigg\rangle_{\omega,\,h} = \bigg\langle\!\!\!\bigg\langle\Delta'_{\omega,\,h}\Omega,\,\Omega\bigg\rangle\!\!\!\bigg\rangle_{\omega,\,h} + \int\limits_X\bigg\langle[i\Theta_h(L)\wedge\cdot\,,\,\Lambda_\omega]\,\Omega,\,\Omega\bigg\rangle_{\omega,\,h}\,dV_\omega,\end{eqnarray*} where the first equality follows from $\langle\langle\Delta''_{\omega,\,h}\Omega,\,\Omega\rangle\rangle_{\omega,\,h} = ||\bar\partial\Omega||^2_{\omega,\,h} + ||\bar\partial^\star_{\omega,\,h}\Omega||^2_{\omega,\,h}$ and from $\bar\partial\Omega = 0$ (by definition of an $s$-symplectic form) and $\bar\partial^\star_{\omega,\,h}\Omega=0$ (for bidegree reasons).

Now, $\langle\langle\Delta'_{\omega,\,h}\Omega,\,\Omega\rangle\rangle_{\omega,\,h} = ||D'_h\Omega||^2_{\omega,\,h} + ||(D'_h)^\star_{\omega,\,h}\Omega||^2_{\omega,\,h}\geq 0$ and, fixing an arbitrary point $x\in X$ and choosing local coordinates $z_1,\dots , z_n$ about $x$ such that (\ref{eqn:coordinates_diagonalisation}) holds, the general identity (\ref{eqn:curvature-operator_eigenvalues}) implies: \begin{eqnarray*}\bigg\langle[i\Theta_h(L)\wedge\cdot\,,\,\Lambda_\omega]\,\Omega,\,\Omega\bigg\rangle_{\omega,\,h} = -\sum\limits_{|J|=s}\lambda_{C_J}\,|\Omega_J^{(\alpha)}|^2\geq 0  \hspace{7ex}\mbox{at}\hspace{1ex} x,\end{eqnarray*} where we put $\Omega_{|U_\alpha} = \sum\limits_{|J| = s}\Omega_J^{(\alpha)}\,dz_J\otimes e_\alpha$ on an open neighbourhood $U_\alpha$ of $x$ on which $L$ is trivialised by $e_\alpha$, $C_J$ is the multi-index that is complementary to $J$ in $\{1,\dots , n\}$, $\lambda_{C_J}$ is the sum of the eigenvalues $\lambda_j$ of $i\Theta_h(L)$ with respect to $\omega$ corresponding to $j\in C_J$ and the last inequality holds if $-\lambda_{C_J}(x)\geq 0$ for every multi-index $J$ of length $|J|=s$.

Since $|C_J| = s$ whenever $|J|=s$ and since $i\Theta_{h^{-1}}(L^{-1}) = -i\Theta_h(L)$, the condition $-\lambda_{C_J}(x)\geq 0$ for every $J$ with $|J|=s$ is equivalent to the sum of any $s$ eigenvalues of the curvature form $i\Theta_{h^{-1}}(L^{-1})$ with respect to $\omega$ being non-negative at $x$. This is further equivalent to the condition $i\Theta_{h^{-1}}(L^{-1})\geq_{s,\,\omega}0$ at $x$.

We conclude that if $i\Theta_{h^{-1}}(L^{-1})\geq_{s,\,\omega}0$ everywhere on $X$, we must have (among other things) $D'_h\Omega = 0$. This last condition is equivalent to \begin{eqnarray}\label{eqn:Omega_divisibility}\partial\Omega_\alpha = \partial\varphi_\alpha\wedge\Omega_\alpha   \hspace{5ex} \mbox{on}\hspace{1ex} U_\alpha \hspace{1ex} \mbox{for every}\hspace{1ex} \alpha,\end{eqnarray} where $\Omega = (\Omega_\alpha)_\alpha$ (i.e. $\Omega_{|U_\alpha} = \Omega_\alpha\otimes e_\alpha$ on $U_\alpha$ for every $\alpha$) and the $\varphi_\alpha$'s are the local weights (each of them defined on the corresponding trivialising open subset $U_\alpha\subset X$) defining the fibre metric $h$ of $L$ (and satisfying the gluing condition (\ref{eqn:gluing_L-h})).

\vspace{1ex}

$\bullet$ If the Hermitian fibre metric $h$ on $L$ discussed above is singular, property (\ref{eqn:Omega_divisibility}) still holds. This can be seen by locally regularising $h$ using Theorem 1.6 in [DP25a] and then using the construction described in $\S$\ref{subsection:Laplacians_singular} as we did in the proof of Theorem \ref{The:p-contact_n-p-pos}. 

\vspace{1ex}

$\bullet$ We shall now see that (\ref{eqn:Omega_divisibility}) leads to a contradiction of the $s$-symplectic assumption on $\Omega$ if $\rk{\cal F}_\Omega>0$. The argument is very similar to the end of the proof of the main result in [Dem02].

The sheaf ${\cal F}_\Omega$ in the statement is a coherent subsheaf of the ${\cal O}_X$-module ${\cal O}(T^{1,\,0}X)$. Property (\ref{eqn:Omega_divisibility}) implies that ${\cal F}_\Omega$ is integrable (i.e. involutive for the Lie bracket $[\,\cdot\,,\,\cdot\,]$ of $T^{1,\,0}X$) thanks to the Cartan formula (third equality below): \begin{eqnarray*}0 & = & (\partial\varphi_\alpha\wedge\Omega_\alpha)(\xi_0,\xi_1,\dots , \xi_s) = (\partial\Omega_\alpha)(\xi_0,\xi_1,\dots , \xi_s) = (d\Omega_\alpha)(\xi_0,\xi_1,\dots , \xi_s)\\
  & = & \sum\limits_{j=0}^s(-1)^j\,\xi_j\cdot\Omega_\alpha(\xi_0, \dots , \widehat{\xi_j}, \dots , \xi_s) + \sum\limits_{0\leq j<k\leq s}(-1)^{j+k}\,\Omega_\alpha\bigg([\xi_j,\,\xi_k],\, \xi_0, \dots , \widehat{\xi_j}, \dots , \widehat{\xi_k}, \dots , \xi_s\bigg)\end{eqnarray*} where these equalities hold for all $(1,\,0)$-vector fields $\xi_0,\xi_1,\dots , \xi_s$, two among which, say $\xi_l$ and $\xi_t$ for some $l<t$, are (local) sections of ${\cal F}_\Omega$. Indeed, we have:

\vspace{1ex}

$\bullet$ the first equality on the first line holds because at most one of $\xi_l$ and $\xi_t$ is acted upon by $\partial\varphi_\alpha$, hence at least one of them contracts $\Omega_\alpha$ (to $0$);

\vspace{1ex}

$\bullet$ all the terms in the first sum on the last line of the above sequence of equalities vanish;

\vspace{1ex}

$\bullet$ all the terms in the second sum vanish, except possibly the term \begin{eqnarray*}(-1)^{l+t}\,\Omega_\alpha\bigg([\xi_l,\,\xi_t],\, \xi_0, \dots , \widehat{\xi_l}, \dots , \widehat{\xi_t}, \dots , \xi_s\bigg).\end{eqnarray*} It follows that this last term must vanish as well for all $(1,\,0)$-vector fields $\xi_0, \dots , \widehat{\xi_l}, \dots , \widehat{\xi_t}, \dots , \xi_s$. This means that $[\xi_l,\,\xi_t]\lrcorner\Omega_\alpha  = 0$, which amounts to $[\xi_l,\,\xi_t]$ being a (local) section of ${\cal F}_\Omega$, for all pairs $\xi_l,\,\xi_t$ of (local) sections of ${\cal F}_\Omega$. This proves that ${\cal F}_\Omega$ is integrable, so it defines a (possibly singular) foliation on $X$.

Now, suppose that $\rk{\cal F}_\Omega>0$. If $x_0\in U_\alpha$ is a smooth point on a leaf of the foliation ${\cal F}_\Omega$, there exist local holomorphic coordinates $z_1,\dots , z_n$ on $U_\alpha$ (shrink $U_\alpha$ about $x_0$ if necessary) such that the leaves of ${\cal F}_\Omega$ are given in $U_\alpha$ by equations \begin{eqnarray*}z_1 = c_1,\dots , z_r=c_r,\end{eqnarray*} where the $c_i$'s are constants and $r: = n - \rk{\cal F}_\Omega <n$. Then, $({\cal F}_\Omega)_{|U_\alpha}$ is generated by $\partial/\partial z_{r+1},\dots , \partial/\partial z_n$ and $\Omega_\alpha$ depends only on $dz_1,\dots , dz_r$. This last fact makes the $s$-symplectic condition \begin{eqnarray*}\Omega_\alpha\wedge\Omega_\alpha \neq 0 \hspace{5ex}\mbox{everywhere on}\hspace{1ex} U_\alpha\end{eqnarray*} impossible for bidegree reasons: $\Omega_\alpha\wedge\Omega_\alpha$ would be a non-vanishing $(n,\,0)$-form depending on only $r$ $dz_j$'s with $r<n=2s$.

  We conclude that $\rk{\cal F}_\Omega = 0$. \hfill $\Box$

\section{The sheaves ${\cal F}_\Gamma$ and ${\cal G}_{\Gamma,\,h}$ associated with a holomorphic $p$-contact bundle}\label{section:F-G_h} These are the analogues in the presence of a $p$-contact line bundle of the sheaves defined in [KPU25] in the scalar setting. While the sheaf ${\cal F}_\Gamma$ continues to be canonically associated with the holomorphic $p$-contact structure $\Gamma$, the sheaf ${\cal G}_{\Gamma,\,h}$ also depends on a fibre metric $h$ on $L$ since $\partial\Gamma$ does not make global sense on $X$.

\begin{Def}\label{Def:F-G_h} Let $(X,\,L,\,\Gamma)$ be a {\bf holomorphic $p$-contact triple} with $\mbox{dim}_\C X = n = 2p+1 = 4l+3$.

\vspace{1ex}

(i)\, We let ${\cal F}_\Gamma$ be the sheaf of germs of $(1,\,0)$-vector fields $\xi$ such that $\xi\lrcorner\Gamma = 0$.

\vspace{1ex}

(ii)\, For every $C^\infty$ Hermitian fibre metric $h$ on $L$, we let ${\cal G}_{\Gamma,\,h}$ be the sheaf of germs of $(1,\,0)$-vector fields $\xi$ such that $\xi\lrcorner D_h'\Gamma = 0$, where $D_h = D_h' + \bar\partial$ is the Chern connection of $(L,\,h)$.

\end{Def}  

The sheaves ${\cal F}_\Gamma$ and ${\cal G}_{\Gamma,\,h}$ are coherent, torsion-free subsheaves of the holomorphic tangent sheaf ${\cal O}(T^{1,\,0}X)$ of $X$ since they are the respective kernels of the holomorphic vector bundle morphisms: \begin{eqnarray*}T^{1,\,0}X\longrightarrow \Lambda^{p-1,\,0}T^\star X\otimes L, & \hspace{5ex} & \xi\longmapsto\xi\lrcorner\Gamma,\\
  T^{1,\,0}X\longrightarrow \Lambda^{p,\,0}T^\star X\otimes L, & \hspace{5ex} & \xi\longmapsto\xi\lrcorner D_h'\Gamma.\end{eqnarray*}

\begin{Prop}\label{Prop:direct-sum_integrability} Let $(X,\,L,\,\Gamma)$ be a {\bf holomorphic $p$-contact triple} with $\mbox{dim}_\C X = n = 2p+1 = 4l+3$. For every $C^\infty$ Hermitian fibre metric $h$ on $L$, the following statements hold.

\vspace{1ex}  

(i)\, The sum \begin{eqnarray*}{\cal F}_\Gamma\oplus {\cal G}_{\Gamma,\,h}\subset{\cal O}(T^{1,\,0}X)\end{eqnarray*} is {\bf direct}.

\vspace{1ex}  

(ii)\, The subsheaf ${\cal G}_{\Gamma,\,h}$ of ${\cal O}(T^{1,\,0}X)$ is {\bf integrable} in the sense that $[{\cal G}_{\Gamma,\,h},\,{\cal G}_{\Gamma,\,h}]\subset{\cal G}_{\Gamma,\,h}$, where $[\,\cdot\,,\,\cdot\,]$ is the Lie bracket of $T^{1,\,0}X$.

\end{Prop}

\noindent {\it Proof.} (i)\, Suppose $\xi\in({\cal F}_\Gamma)_x\cap({\cal G}_{\Gamma,\,h})_x$ for some $x\in X$. Then, $\xi\lrcorner\Gamma = 0$ and $\xi\lrcorner D_h'\Gamma = 0$ on some neighbourhood $U$ of $x$. This implies that $\xi\lrcorner(\Gamma\wedge D_h'\Gamma) = (\xi\lrcorner\Gamma)\wedge D_h'\Gamma - \Gamma\wedge(\xi\lrcorner D_h'\Gamma) = 0$ on $U$.

Now, recall that $\Gamma\wedge D_h'\Gamma = \Gamma\wedge\partial\Gamma = (\Gamma_\alpha\wedge\partial\Gamma_\alpha)_\alpha\in C^\infty_{n,\,0}(X,\,L^2)$ (see Observation \ref{Obs:Gamma-wedge-sel-Gamma_global}), so the condition $\xi\lrcorner(\Gamma\wedge D_h'\Gamma) = 0$ on $U$ is equivalent to $\xi\lrcorner(\Gamma_\alpha\wedge\partial\Gamma_\alpha) = 0$ on $U\cap U_\alpha$ for every $\alpha$. Since $\Gamma_\alpha\wedge\partial\Gamma_\alpha$ is a {\it non-vanishing} scalar $(n,\,0)$-form on $U_\alpha$, this is further equivalent to $\xi = 0$ on $U\cap U_\alpha$ (for every $\alpha$).

\vspace{1ex}  

(ii)\, Since $D_h(D_h'\Gamma) = \bar\partial D_h'\Gamma = (D_h'\bar\partial + \bar\partial D_h')\,\Gamma = \Theta_h(L)\wedge\Gamma$ is an $L$-valued form of bidegree $(p+1,\,1)$, it vanishes on any $(p+2)$-tuple $(\xi_0,\xi_1,\dots , \xi_{p+1})$ of $(1,\,0)$-vector fields. Hence, we get the first equality below, where the second equality follows from Cartan's formula involving the covariant derivative $\nabla^{(h)}_{\xi_j}$ induced by the connection $D_h$ in the direction of any $(1,\,0)$-vector field $\xi_j$: \begin{eqnarray}\label{eqn:Cartan-application_G-h}0 = D_h(D_h'\Gamma)(\xi_0,\xi_1,\dots , \xi_{p+1}) & = & \sum\limits_{j=0}^{p+1}(-1)^j\,\nabla^{(h)}_{\xi_j}\bigg((D'_h\Gamma)(\xi_0, \dots , \widehat{\xi_j}, \dots , \xi_{p+1})\bigg)\\
\nonumber  & + & \sum\limits_{0\leq j<k\leq p+1}(-1)^{j+k}\,(D'_h\Gamma)\bigg([\xi_j,\,\xi_k],\, \xi_0, \dots , \widehat{\xi_j}, \dots , \widehat{\xi_k}, \dots , \xi_{p+1}\bigg)\end{eqnarray} for all $(1,\,0)$-vector fields $\xi_0,\xi_1,\dots , \xi_{p+1}$.

If two among $\xi_0,\xi_1,\dots , \xi_{p+1}$, say $\xi_l$ and $\xi_s$ for some $l<s$, are (local) sections of ${\cal G}_{\Gamma,\,h}$, then:

\vspace{1ex}

$\bullet$ all the terms in the first sum on the r.h.s. of (\ref{eqn:Cartan-application_G-h}) vanish;

\vspace{1ex}

$\bullet$ all the terms in the second sum on the r.h.s. of (\ref{eqn:Cartan-application_G-h}) vanish, except possibly the term \begin{eqnarray*}(-1)^{l+s}\,(D'_h\Gamma)\bigg([\xi_l,\,\xi_s],\, \xi_0, \dots , \widehat{\xi_l}, \dots , \widehat{\xi_s}, \dots , \xi_{p+1}\bigg).\end{eqnarray*}

\vspace{1ex}

It follows that this last term must vanish as well for all $(1,\,0)$-vector fields $\xi_0, \dots , \widehat{\xi_l}, \dots , \widehat{\xi_s}, \dots , \xi_{p+1}$. This means that $[\xi_l,\,\xi_s]\lrcorner D'_h\Gamma = 0$, which amounts to $[\xi_l,\,\xi_s]$ being a (local) section of ${\cal G}_{\Gamma,\,h}$, for all pairs $\xi_l,\,\xi_s$ of (local) sections of ${\cal G}_{\Gamma,\,h}$. Thus, ${\cal G}_{\Gamma,\,h}$ is integrable.  \hfill $\Box$

\vspace{2ex}

Part (ii) of Proposition \ref{Prop:direct-sum_integrability} can be reworded into the statement that ${\cal G}_{\Gamma,\,h}$ defines a (possibly singular) foliation on $X$.

\section{Examples of (locally holomorphic) $p$-contact manifolds}\label{section:examples}

We will point out examples in several well-known classes of compact complex spin manifolds.

\subsection{Complex projective spaces}\label{subsection:projective-spaces}

For any positive integer $n$, let $\Proj^n = \C\Proj^n$ be the $n$-dimensional complex projective space. Its canonical bundle is $K_{\Proj^n} = {\cal O}(-n-1)$, so we deduce from the classical Proposition \ref{Prop:spin_existence_standard} that $\Proj^n$ admits a spin structure if and only if $n$ is odd. This, the dimension constraints and Corollaries \ref{Cor:p-contact_spin} and \ref{Cor:s-symplectic_spin} yield the following statements:

\vspace{1ex}

$\bullet$ No complex projective space $\Proj^n$ is a (locally holomorphic) $s$-symplectic manifold.

\vspace{1ex}

$\bullet$ For a projective space $\Proj^n$ to be a (locally holomorphic) $p$-contact manifold it is necessary that $n\equiv 3$ mod $4$.

\vspace{1ex}

We will show that this necessary condition is also sufficient. Let $n=2p+1$ with $p$ odd. Any local holomorphic $p$-contact structure $\Gamma$ on $\Proj^n$ will take values in a line bundle $L={\cal O}(k)$ (since the ${\cal O}(k)$'s with $k\in\Z$ are the only holomorphic line bundles on $\Proj^n$) such that $L^2 = {\cal O}(2k)\simeq -K_{\Proj^n} = {\cal O}(n+1) = {\cal O}(2p+2)$. This forces $k=p+1$. Thus, the first thing we have to check is whether the vector space $H^{p,\,0}(\Proj^n,\,{\cal O}(p+1))$ is non-zero.

This is, indeed, known to be true (see, e.g. [Dem97, chapter VII, $\S10$]). Recall that if $V$ is a $\C$-vector space with $\mbox{dim}_\C V = n+1$ on which coordinates $z=(z_0,\, z_1,\dots , z_n)$ have been fixed, one denotes by \begin{eqnarray*}\xi=\sum\limits_{j=0}^n z_j\,\frac{\partial}{\partial z_j}\end{eqnarray*} the {\it Euler vector field} and gets (cf. Theorem 10.6 in [Dem97, chapter VII]): \begin{eqnarray}\label{eqn:H_p0_O(k)}\nonumber H^{p,\,0}(\Proj(V),\,{\cal O}(k)) & \simeq & Z^{p,\,k}(V^\star):=\bigg\{\alpha\in \Lambda^p V^\star\otimes S^{k-p}V^\star\,\mid\,\xi\lrcorner\alpha = 0\bigg\}  \hspace{3ex} \mbox{for all}\hspace{1ex} k\geq p\geq 0 \\
  H^{p,\,0}(\Proj(V),\,{\cal O}(k)) & = & 0 \hspace{25ex} \mbox{whenever}\hspace{1ex} k\leq p \hspace{1ex}\mbox{and}\hspace{1ex} (k,\,p)\neq(0,\,0).\end{eqnarray} There is, of course, a natural identification: \begin{eqnarray*}\Lambda^p V^\star\otimes S^{k-p}V^\star\simeq\bigg\{\alpha=\sum\limits_{|I|=p}\alpha_I(z)\,dz_I\,\mid\,\alpha_I \hspace{2ex}\mbox{are homogeneous polynomials of degree}\hspace{1ex} k-p\bigg\}.\end{eqnarray*}

The case we are interested in (i.e. $k=p+1$) lies on the first row of (\ref{eqn:H_p0_O(k)}), so we get \begin{eqnarray*}H^{p,\,0}(\Proj^n,\,{\cal O}(p+1))\simeq Z^{p,\,p+1}\bigg((\C^{n+1})^\star\bigg)\neq 0.\end{eqnarray*}

We can now go further by running an explicit construction and prove the following

\begin{Prop}\label{Prop:proj-space_examples} For every positive integer $n\equiv 3$ mod $4$, the complex projective space $\Proj^n$ is a (locally holomorphic) {\bf $p$-contact manifold}.

\end{Prop}  

\noindent {\it Proof.} Let $n=2p+1$ with $p$ {\it odd}. Let $(x_0,\, x_1,\dots , x_n)$ be the Euclidean holomorphic coordinates on $\C^{n+1}$ and, for $\alpha=0,1,\dots , n$, let $U_\alpha = \{[x_0:\, x_1:\dots : x_n]\in\Proj^n\,\mid\,x_\alpha\neq 0\}\simeq\C^n$ be the associated open coordinate patches of $\Proj^n$. The induced holomorphic coordinates on each $U_\alpha$ are: \begin{eqnarray*}z_0^{(\alpha)} = \frac{x_0}{x_\alpha},\dots , z_{\alpha-1}^{(\alpha)} = \frac{x_{\alpha-1}}{x_\alpha}, \,z_{\alpha+1}^{(\alpha)} = \frac{x_{\alpha+1}}{x_\alpha},\dots , z_n^{(\alpha)} = \frac{x_n}{x_\alpha}.\end{eqnarray*}      

We will explicitly construct a global holomorphic form $\Gamma = (\Gamma_\alpha)_{0\leq\alpha\leq n}\in H^{p,\,0}(\Proj^n,\,{\cal O}(p+1))$ such that $\Gamma\wedge\partial\Gamma\neq 0$ everywhere on $\Proj^n$ (i.e. $\Gamma_\alpha\wedge\partial\Gamma_\alpha\neq 0$ everywhere on $U_\alpha$ for every $\alpha$). To this end, we will first construct $\Gamma_0$ on $U_0$ and will then compute $\Gamma_\alpha$ for every $\alpha\neq 0$ such that $\Gamma_0 = g_{0\alpha}^{(p+1)}\,\Gamma_\alpha$ on $U_0\cap U_\alpha$, where the $g_{\alpha\beta}^{(p+1)}$'s are the transition functions of the holomorphic line bundle ${\cal O}(p+1)$ relative to its standard trivialisations over the $U_\alpha$'s.

The first step is the computation of the $g_{0\alpha}^{(p+1)}$'s and of the way the products of $p$ $dz_j^{(0)}$'s transform when passing from $U_0$ to $U_\alpha$.

\vspace{1ex}

$\bullet$ We fix an arbitrary $\alpha\in\{1,\dots, n\}$. Since ${\cal O}(p+1)_{|U_0}$ is holomorphically trivialised by the section $x_0^{p+1}$ and ${\cal O}(p+1)_{|U_\alpha}$ is holomorphically trivialised by the section $x_\alpha^{p+1}$ and \begin{eqnarray*}x_\alpha^{p+1} = \bigg(\frac{x_\alpha}{x_0}\bigg)^{p+1}\,x_0^{p+1} = \bigg(\frac{1}{z_0^{(\alpha)}}\bigg)^{p+1}\,x_0^{p+1}  \hspace{5ex} \mbox{on}\hspace{2ex} U_0\cap U_\alpha,\end{eqnarray*} we deduce that $g_{0\alpha}^{(p+1)} = (1/z_0^{(\alpha)})^{p+1}$ on $U_0\cap U_\alpha$.

Now, let $J:=(j_1,\dots, j_p)$ with $1\leq j_1<\dots <j_p\leq n$. If $j_k=\alpha$ for some $k$, we get: \begin{eqnarray*}z_{j_1}^{(0)} = \frac{z_{j_1}^{(\alpha)}}{z_0^{(\alpha)}}, \dots , z_{j_{k-1}}^{(0)} = \frac{z_{j_{k-1}}^{(\alpha)}}{z_0^{(\alpha)}},\, z_{j_k}^{(0)} = \frac{1}{z_0^{(\alpha)}},\,z_{j_{k+1}}^{(0)} = \frac{z_{j_{k+1}}^{(\alpha)}}{z_0^{(\alpha)}},\dots, z_{j_p}^{(0)} = \frac{z_{j_p}^{(\alpha)}}{z_0^{(\alpha)}},\end{eqnarray*} hence \begin{eqnarray*}dz_{j_1}^{(0)} = \frac{1}{z_0^{(\alpha)}}\, dz_{j_1}^{(\alpha)} - \frac{z_{j_1}^{(\alpha)}}{(z_0^{(\alpha)})^2}\, dz_0^{(\alpha)}, \dots , dz_{j_k}^{(0)} = -\frac{1}{(z_0^{(\alpha)})^2}\,dz_0^{(\alpha)}, \dots , dz_{j_p}^{(0)} = \frac{1}{z_0^{(\alpha)}}\, dz_{j_p}^{(\alpha)} - \frac{z_{j_p}^{(\alpha)}}{(z_0^{(\alpha)})^2}\, dz_0^{(\alpha)}\end{eqnarray*} on $U_0\cap U_\alpha$.

These relations imply that, when $\alpha\in\{j_1,\dots , j_p\}$, say $\alpha=j_k$, we get: \begin{eqnarray}\label{eqn:dz_J^0_index-inside}dz_{j_1}^{(0)}\wedge\dots\wedge dz_{j_p}^{(0)} = \frac{1}{(z_0^{(\alpha)})^{p+1}}\,(-1)^k\,dz_0^{(\alpha)}\wedge dz_{j_1}^{(\alpha)}\wedge\dots\wedge\widehat{dz_{j_k}^{(\alpha)}}\wedge\dots\wedge dz_{j_p}^{(\alpha)}\end{eqnarray} on $U_0\cap U_\alpha$ because the factor $dz_{j_k}^{(0)}$ on the left is a function multiple of $dz_0^{(\alpha)}$ and this selects only the first term, a function multiple of $dz_{j_l}^{(\alpha)}$, from each expression for $dz_{j_l}^{(0)}$ with $l\neq k$.

When $\alpha\notin\{j_1,\dots , j_p\}$, we get: \begin{eqnarray*} dz_{j_1}^{(0)}\wedge\dots\wedge dz_{j_p}^{(0)} = \frac{1}{(z_0^{(\alpha)})^p}\,dz_{j_1}^{(\alpha)}\wedge\dots\wedge dz_{j_p}^{(\alpha)} & - & \frac{z_{j_1}^{(\alpha)}}{(z_0^{(\alpha)})^{p+1}}\,dz_0^{(\alpha)}\wedge dz_{j_2}^{(\alpha)}\wedge\dots\wedge dz_{j_p}^{(\alpha)}\\
\nonumber & - & \dots - \frac{z_{j_p}^{(\alpha)}}{(z_0^{(\alpha)})^{p+1}}\, dz_{j_1}^{(\alpha)}\wedge\dots\wedge dz_{j_{p-1}}^{(\alpha)}\wedge dz_0^{(\alpha)}\end{eqnarray*} on $U_0\cap U_\alpha$, where the first term on the right is obtained by selecting the first term, a function multiple of $dz_{j_l}^{(\alpha)}$, from each expression for $dz_{j_l}^{(0)}$, while all the other terms on the right are obtained by selecting the second term, a function multiple of $dz_0^{(\alpha)}$, from the expression of one $dz_{j_l}^{(0)}$ and the first term, a function multiple of $dz_{j_r}^{(\alpha)}$, from the expression for $dz_{j_r}^{(0)}$ for each $r\neq l$. Thus, for $\alpha\notin\{j_1,\dots , j_p\}$, we get: \begin{eqnarray}\label{eqn:dz_J^0_index-outside} dz_{j_1}^{(0)}\wedge\dots\wedge dz_{j_p}^{(0)} & = & \frac{1}{(z_0^{(\alpha)})^{p+1}}\,\bigg((z_0^{(\alpha)}\,dz_{j_1}^{(\alpha)}\wedge\dots\wedge dz_{j_p}^{(\alpha)} + \sum\limits_{l=1}^p(-1)^l\,z_{j_l}^{(\alpha)}\,\widehat{dz_{j_l}^{(\alpha)}}\bigg)\end{eqnarray} on $U_0\cap U_\alpha$, where \begin{eqnarray*}\widehat{dz_{j_l}^{(\alpha)}}:= dz_0^{(\alpha)}\wedge dz_{j_1}^{(\alpha)}\wedge\dots\wedge\widehat{dz_{j_l}^{(\alpha)}}\wedge\dots\wedge dz_{j_p}^{(\alpha)}, \hspace{5ex} l=1,\dots , p.\end{eqnarray*}

\vspace{1ex}

$\bullet$ Let us now choose $\Gamma_0:=a_0\,dz_1^{(0)}\wedge\dots\wedge dz_p^{(0)} + b_0\,dz_{p+1}^{(0)}\wedge\dots\wedge dz_{n-1}^{(0)}\in C^\infty_{p,\,0}(U_0,\,\C)$ with holomorphic functions $a_0,\, b_0:U_0\longrightarrow\C$. Then: \begin{eqnarray*}\Gamma_0\wedge\partial\Gamma_0 & = & \bigg(b_0\,\frac{\partial a_0}{\partial z_n^{(0)}} + (-1)^p\,a_0\,\frac{\partial b_0}{\partial z_n^{(0)}}\bigg)\,dz_1^{(0)}\wedge\dots\wedge dz_n^{(0)} \\
  & = & b_0^2\,\frac{\partial}{\partial z_n^{(0)}}\bigg(\frac{a_0}{b_0}\bigg)\,dz_1^{(0)}\wedge\dots\wedge dz_n^{(0)}     \hspace{5ex} \mbox{on}\hspace{2ex} U_0,\end{eqnarray*} the last equality holding thanks to $p$ being odd if $b_0$ is non-vanishing. Since we want $\Gamma_0\wedge\partial\Gamma_0$ to be non-vanishing on $U_0$, a natural choice is $b_0\equiv 1$ and $a_0 = z_n^{(0)}$, which yields: \begin{eqnarray}\label{eqn:Gamma_0}\Gamma_0:=z_n^{(0)}\,dz_1^{(0)}\wedge\dots\wedge dz_p^{(0)} + dz_{p+1}^{(0)}\wedge\dots\wedge dz_{n-1}^{(0)}  \hspace{5ex} \mbox{on}\hspace{2ex} U_0,\end{eqnarray} hence \begin{eqnarray*}\Gamma_0\wedge\partial\Gamma_0 = dz_1^{(0)}\wedge\dots\wedge dz_n^{(0)} \neq 0  \hspace{5ex} \mbox{everywhere on}\hspace{2ex} U_0.\end{eqnarray*}

\vspace{1ex}

$\bullet$ Let $\alpha\in\{1,\dots , p\}$. From (\ref{eqn:dz_J^0_index-inside}) and (\ref{eqn:dz_J^0_index-outside}), we get on $U_0\cap U_\alpha$: \begin{eqnarray*}\Gamma_0 =  \frac{1}{(z_0^{(\alpha)})^{p+1}}\,\bigg(& & (-1)^\alpha\,\frac{z_n^{(\alpha)}}{z_0^{(\alpha)}}\, dz_0^{(\alpha)}\wedge dz_1^{(\alpha)}\wedge\dots\wedge\widehat{dz_\alpha^{(\alpha)}}\wedge\dots\wedge dz_p^{(\alpha)} \\
  & + & z_0^{(\alpha)}\,dz_{p+1}^{(\alpha)}\wedge\dots\wedge dz_{n-1}^{(\alpha)} - z_{p+1}^{(\alpha)}\,dz_0^{(\alpha)}\wedge dz_{p+2}^{(\alpha)}\wedge\dots\wedge dz_{n-1}^{(\alpha)} \\
  & + & z_{p+2}^{(\alpha)}\,dz_0^{(\alpha)}\wedge dz_{p+1}^{(\alpha)}\wedge dz_{p+3}^{(\alpha)}\wedge\dots\wedge dz_{n-1}^{(\alpha)}  \\
  & \vdots & \\
  & + &  (-1)^p\, z_{n-1}^{(\alpha)}\,dz_0^{(\alpha)}\wedge dz_{p+1}^{(\alpha)}\wedge\dots\wedge dz_{n-2}^{(\alpha)}\bigg):= \frac{1}{(z_0^{(\alpha)})^{p+1}}\,\Gamma_\alpha,\end{eqnarray*} where the last equality constitutes the definition of a $\C$-valued holomorphic $(p,\,0)$-form $\Gamma_\alpha$ on $U_\alpha$ that links to the given $\Gamma_0$ on $U_0\cap U_\alpha$ via the transition function $g_{0\alpha}^{(p+1)}$ of ${\cal O}(p+1)$. This $\Gamma_\alpha$ is forced on us by $\Gamma_0$ and the coordinate changes computed above.

We now compute and get: \begin{eqnarray*}\Gamma_\alpha\wedge\partial\Gamma_\alpha & = & (-1)^\alpha\,\frac{\partial\bigg(\frac{z_n^{(\alpha)}}{z_0^{(\alpha)}}\bigg)}{\partial z_n^{(\alpha)}}\,z_0^{(\alpha)}\,dz_n^{(\alpha)}\wedge dz_0^{(\alpha)}\wedge\dots\wedge\widehat{dz_\alpha^{(\alpha)}}\wedge\dots\wedge dz_{n-1}^{(\alpha)} \\
  & = & (-1)^\alpha\,dz_n^{(\alpha)}\wedge dz_0^{(\alpha)}\wedge\dots\wedge\widehat{dz_\alpha^{(\alpha)}}\wedge\dots\wedge dz_{n-1}^{(\alpha)} \\
  & = & (-1)^{n+\alpha-1}\,dz_0^{(\alpha)}\wedge\dots\wedge\widehat{dz_\alpha^{(\alpha)}}\wedge\dots\wedge dz_n^{(\alpha)} \neq 0\end{eqnarray*} everywhere on $U_\alpha$ for every $\alpha\in\{1,\dots , p\}$.

\vspace{1ex}

$\bullet$ Let $\alpha\in\{p+1,\dots , n-1\}$. From (\ref{eqn:dz_J^0_index-inside}) and (\ref{eqn:dz_J^0_index-outside}), we get on $U_0\cap U_\alpha$: \begin{eqnarray*}\Gamma_0 =  \frac{1}{(z_0^{(\alpha)})^{p+1}}\,\bigg(& & \frac{z_n^{(\alpha)}}{z_0^{(\alpha)}}\,z_0^{(\alpha)}\,dz_1^{(\alpha)}\wedge\dots\wedge dz_p^{(\alpha)} - \frac{z_n^{(\alpha)}}{z_0^{(\alpha)}}\,z_1^{(\alpha)}\,dz_0^{(\alpha)}\wedge dz_2^{(\alpha)}\wedge\dots\wedge dz_p^{(\alpha)} \\
  & + & \frac{z_n^{(\alpha)}}{z_0^{(\alpha)}}\,z_2^{(\alpha)}\,dz_0^{(\alpha)}\wedge\dots\wedge\widehat{dz_2^{(\alpha)}}\wedge\dots\wedge dz_p^{(\alpha)} \\
  & \vdots & \\
  & + & (-1)^p\,\frac{z_n^{(\alpha)}}{z_0^{(\alpha)}}\,z_p^{(\alpha)}\,dz_0^{(\alpha)}\wedge\dots\wedge dz_{p-1}^{(\alpha)} \\
  & + & (-1)^{\alpha-p}\,dz_0^{(\alpha)}\wedge dz_{p+1}^{(\alpha)}\wedge\dots\wedge\widehat{dz_\alpha^{(\alpha)}}\wedge\dots\wedge dz_{n-1}^{(\alpha)}\bigg):= \frac{1}{(z_0^{(\alpha)})^{p+1}}\,\Gamma_\alpha,\end{eqnarray*} where the last equality constitutes the definition of a $\C$-valued holomorphic $(p,\,0)$-form $\Gamma_\alpha$ on $U_\alpha$ that links to the given $\Gamma_0$ on $U_0\cap U_\alpha$ via the transition function $g_{0\alpha}^{(p+1)}$ of ${\cal O}(p+1)$.

We now compute and get: \begin{eqnarray*}\Gamma_\alpha\wedge\partial\Gamma_\alpha & = & (-1)^{\alpha-p}\,\frac{\partial z_n^{(\alpha)}}{\partial z_n^{(\alpha)}}\,dz_n^{(\alpha)}\wedge dz_1^{(\alpha)}\wedge\dots\wedge dz_p^{(\alpha)}\wedge dz_0^{(\alpha)}\wedge dz_{p+1}^{(\alpha)}\wedge\dots\wedge\widehat{dz_\alpha^{(\alpha)}}\wedge\dots\wedge dz_{n-1}^{(\alpha)}  \\
  & = & (-1)^{n+\alpha-1}\,dz_0^{(\alpha)}\wedge\dots\wedge\widehat{dz_\alpha^{(\alpha)}}\wedge\dots\wedge dz_n^{(\alpha)}\neq 0\end{eqnarray*} everywhere on $U_\alpha$ for every $\alpha\in\{p+1,\dots , n-1\}$.

\vspace{1ex}

$\bullet$ Let $\alpha = n$. From (\ref{eqn:dz_J^0_index-outside}), we get on $U_0\cap U_n$: \begin{eqnarray*}\Gamma_0 =  \frac{1}{(z_0^{(n)})^{p+1}}\,\bigg(& & \frac{1}{z_0^{(n)}}\,z_0^{(n)}\,dz_1^{(n)}\wedge\dots\wedge dz_p^{(n)} - \frac{1}{z_0^{(n)}}\,z_1^{(n)}\,dz_0^{(n)}\wedge dz_2^{(n)}\wedge\dots\wedge dz_p^{(n)} \\
 & + & \frac{1}{z_0^{(n)}}\,z_2^{(n)}\,dz_0^{(n)}\wedge\dots\wedge\widehat{dz_2^{(n)}}\wedge\dots\wedge dz_p^{(n)} \\
  & \vdots & \\
  & + & (-1)^p\,\frac{1}{z_0^{(n)}}\,z_p^{(n)}\,dz_0^{(n)}\wedge\dots\wedge dz_{p-1}^{(n)} \\
  & + & z_0^{(n)}\,dz_{p+1}^{(n)}\wedge\dots\wedge dz_{n-1}^{(n)} -  z_{p+1}^{(n)}\,dz_0^{(n)}\wedge dz_{p+2}^{(n)}\wedge\dots\wedge dz_{n-1}^{(n)} \\
  & \vdots & \\
  & + & (-1)^p\,z_{n-1}^{(n)}\,dz_0^{(n)}\wedge dz_{p+1}^{(n)}\wedge\dots\wedge dz_{n-2}^{(n)}\bigg):= \frac{1}{(z_0^{(n)})^{p+1}}\,\Gamma_n,\end{eqnarray*} where the last equality constitutes the definition of a $\C$-valued holomorphic $(p,\,0)$-form $\Gamma_n$ on $U_n$ that links to the given $\Gamma_0$ on $U_0\cap U_n$ via the transition function $g_{0n}^{(p+1)}$ of ${\cal O}(p+1)$. 

We now compute and get: \begin{eqnarray*}\Gamma_n\wedge\partial\Gamma_n & = & -\frac{z_0^{(n)}}{z_0^{(n)}}\,dz_1^{(n)}\wedge dz_0^{(n)}\wedge dz_2^{(n)}\wedge\dots\wedge dz_{n-1}^{(n)} + \frac{z_0^{(n)}}{z_0^{(n)}}\,dz_2^{(n)}\wedge dz_0^{(n)}\wedge\dots\wedge\widehat{dz_2^{(n)}}\wedge\dots\wedge dz_{n-1}^{(n)}   \\
  & \vdots & \\
  & + & \frac{z_0^{(n)}}{z_0^{(n)}}\,(-1)^p\,dz_p^{(n)}\wedge dz_0^{(n)}\wedge\dots\wedge\widehat{dz_p^{(n)}}\wedge\dots\wedge dz_{n-1}^{(n)} \\
  & = & p\,dz_0^{(n)}\wedge\dots\wedge dz_{n-1}^{(n)}\neq 0\end{eqnarray*} everywhere on $U_n$.

\vspace{1ex}

$\bullet$ We have thus checked that $\Gamma_\alpha\wedge\partial\Gamma_\alpha \neq 0$ everywhere on $U_\alpha$ for every $\alpha\in\{1,\dots , n\}$. This proves that $\Gamma\wedge\partial\Gamma\neq 0$ everywhere on $\Proj^n$, completing the proof.  \hfill $\Box$

\subsection{Smooth hypersurfaces in complex projective spaces}\label{subsection:hypersurfaces-projective-spaces}

Let $X\subset\Proj^{n+1}$ be a smooth complex hypersurface of degree $d$ in the $(n+1)$-dimensional complex projective space. This means that $X$ is the set of zeroes of a degree-$d$ homogeneous polynomial $P$ on $\C^{n+2}$, identified with a global holomorphic section of the line bundle ${\cal O}(d) = {\cal O}_{\Proj^{n+1}}(d)$, such that $(dP)(z)\neq 0$ for every $z=[z_0:\dots : z_{n+1}]\in\Proj^{n+1}$ for which $P(z)=0$.

It is standard that the anti-canonical bundle of $X$ is given by $-K_X = {\cal O}_{\Proj^{n+1}}(n+2-d)_{|X}$, so $X$ is a spin manifold if and only if $n+2-d$ is even. When $n$ is odd, this is equivalent to $d$ being odd as well.

Let us now suppose that $\mbox{dim}_\C X = n = 2p+1 = 4l+3$, with $p = 2l+1$ for some non-negative integer $l$. As we have seen, whenever $d$ is odd, $X$ admits a spin structure that we can identify with a holomorphic line bundle $L$ on $X$ such that $L^2\simeq -K_X = {\cal O}_{\Proj^{n+1}}(n+2-d)_{|X} = {\cal O}_X(n+2-d)$. The obvious choice is $L={\cal O}_X((n+2-d)/2)$, which means that \begin{eqnarray}\label{eqn:L_hypersurface_Pn+1}L={\cal O}_X(k), \hspace{3ex} \mbox{with}\hspace{2ex} k = p+1 + \frac{1-d}{2}.\end{eqnarray}

An $L$-valued holomorphic $p$-contact structure on $X$, if any, is a form \begin{eqnarray}\label{eqn:Gamma_L_hypersurface_Pn+1}\Gamma\in H^{p,\,0}\bigg(X,\,{\cal O}_X\bigg(p+1 + \frac{1-d}{2}\bigg)\bigg) \hspace{3ex}\mbox{such that}\hspace{2ex} \Gamma\wedge\partial\Gamma\neq 0 \hspace{3ex}\mbox{everywhere on}\hspace{2ex} X.\end{eqnarray}

Meanwhile, we know from Corollary \ref{Cor:Demailly_K_not-psef} (a consequence of the proof of the main result of [Dem02]) that if an $L$-valued holomorphic $p$-contact structure $\Gamma$ exists, $K_X = {\cal O}_X(d-n-2)$ cannot be pseudo-effective. This forces $d<n+2$. Since $d$ is odd and $n+1$ is even, this means that $d\leq n$. The case $d=1$ corresponds to the smooth complex hypersurface $X\subset\Proj^{n+1}$ of degree $d$ being a copy of $\Proj^n$ which, thanks to Proposition \ref{Prop:proj-space_examples}, is a (locally holomorphic) {\bf $p$-contact manifold} when $n=2p+1 = 4l+3$.

Unfortunately, as we shall now see, the case $d=1$ is the only case where a hypersurface $X\subset\Proj^{n+1}$ has a $p$-contact structure.

\begin{Prop}\label{Prop:hypersurfaces-projective-spaces} Let $n=2p+1 = 4l+3$ and let $X\subset\Proj^{n+1}$ be a smooth complex hypersurface of odd degree $d$ such that $d\leq n$. Set $k = (n+2-d)/2 = p+1 + (1-d)/2.$

  If $d\geq 3$, the following cohomology vanishing holds: \begin{eqnarray}\label{eqn:hypersurfaces-projective-spaces_coh-vanishing}H^{p,\,0}(X,\,{\cal O}_X(k)) = 0.\end{eqnarray}

  In particular, there is {\bf no holomorphic $p$-contact structure} on $X$ if $\mbox{deg}\,X\geq 3$. In other words, the only smooth hypersurfaces in complex projective space that carry a holomorphic $p$-contact structure are copies of $\Proj^n$ with $n\equiv 3$ mod $4$.

\end{Prop}  

\noindent {\it Proof.}\footnote{The second-named author is grateful to L. Manivel for pointing out to him the main arguments in this proof.} $\bullet$ {\it Step $1$.} Since $\deg X =d$, the normal bundle (a holomorphic line bundle) of $X$ in $\Proj^{n+1}$ is ${\cal N}_{X/\Proj^{n+1}} = {\cal O}_X(d)$, so the {\it normal} short exact sequence of holomorphic vector bundles over $X$ reads:\begin{eqnarray*}0\longrightarrow T^{1,\,0}X \longrightarrow T^{1,\,0}\Proj^{n+1}_{|X}\longrightarrow {\cal O}_X(d) \longrightarrow 0.\end{eqnarray*} The dual short exact sequence is \begin{eqnarray}\label{eqn:dual-normal-sequence_X}0\longrightarrow {\cal O}_X(-d) = \Omega^0_X(-d) \longrightarrow \Lambda^{1,\,0}T^\star\Proj^{n+1}_{|X} = \Omega^1_{\Proj^{n+1}|X}\longrightarrow \Lambda^{1,\,0}T^\star X = \Omega^1_X\longrightarrow 0,\end{eqnarray} where the bundle notation and the sheaf notation are used interchangeably.

Tensoring (\ref{eqn:dual-normal-sequence_X}) by the holomorphic vector bundle $\Omega^{p-1}_X(k)$, we get the short exact sequence of vector bundles over $X$: \begin{eqnarray*}\label{eqn:normal-sequence_X}0\longrightarrow \Omega^{p-1}_X(k-d) \longrightarrow \Omega^p_{\Proj^{n+1}|X}(k)\longrightarrow \Omega^p_X(k)\longrightarrow 0.\end{eqnarray*} From the induced long exact sequence in cohomology, we extract \begin{eqnarray*}H^0(X,\,\Omega^p_{\Proj^{n+1}|X}(k))\longrightarrow H^0(X,\,\Omega^p_X(k)) \longrightarrow H^1(X,\,\Omega^{p-1}_X(k-d)).\end{eqnarray*} This exact sequence shows that if $H^0(X,\,\Omega^p_{\Proj^{n+1}|X}(k)) = 0$, the vector space $H^0(X,\,\Omega^p_X(k))$ injects into $H^1(X,\,\Omega^{p-1}_X(k-d))$.

More generally, tensoring (\ref{eqn:dual-normal-sequence_X}) by the holomorphic vector bundle $\Omega^{p-i-1}_X(k-id)$, for $i\in\{0,\dots, p-1\}$, we get the short exact sequence of vector bundles over $X$: \begin{eqnarray*}\label{eqn:normal-sequence_X}0\longrightarrow \Omega^{p-i-1}_X(k-(i+1)d) \longrightarrow \Omega^{p-i}_{\Proj^{n+1}|X}(k-id)\longrightarrow \Omega^{p-i}_X(k-id)\longrightarrow 0.\end{eqnarray*} From the induced long exact sequence in cohomology, we extract \begin{eqnarray*}H^i(X,\,\Omega^{p-i}_{\Proj^{n+1}|X}(k-id))\longrightarrow H^i(X,\,\Omega^{p-i}_X(k-id)) \longrightarrow H^{i+1}(X,\,\Omega^{p-i-1}_X(k-(i+1)d)).\end{eqnarray*} This exact sequence shows that if $H^i(X,\,\Omega^{p-i}_{\Proj^{n+1}|X}(k-id)) = 0$, the vector space $H^i(X,\,\Omega^{p-i}_X(k-id))$ injects into $H^{i+1}(X,\,\Omega^{p-i-1}_X(k-(i+1)d))$.

In other words, this shows that if \begin{eqnarray}\label{eqn:multiple-vanishings}H^0(X,\,\Omega^p_{\Proj^{n+1}|X}(k)) = \dots = H^i(X,\,\Omega^{p-i}_{\Proj^{n+1}|X}(k-id)) = \dots = H^{p-1}(X,\,\Omega^1_{\Proj^{n+1}|X}(k-(p-1)d)) = 0,\end{eqnarray} we get the sequence of injections:  \begin{eqnarray*}H^0(X,\,\Omega^p_X(k)) & \hookrightarrow & \dots \hookrightarrow H^i(X,\,\Omega^{p-i}_X(k-id))\hookrightarrow  H^{i+1}(X,\,\Omega^{p-i-1}_X(k-(i+1)d))\hookrightarrow \dots \\
  & \hookrightarrow & H^p(X,\,\Omega^0_X(k-pd)) = H^{0,\,p}(X,\,{\cal O}_X(k-pd)),\end{eqnarray*} where ${\cal O}_X(k-pd))$ denotes a line bundle rather than the associated sheaf.

Now, the hypothesis $d\leq n$ implies that $X$ is a Fano manifold since, as we have seen above, $-K_X = {\cal O}_X(n+2-d)$ and $n+2-d>0$. This further implies that $H^{0,\,p}(X,\,{\cal O}_X(k-pd)) = 0$ whatever the sign of $k-pd$. Indeed, if $k-pd < 0$, this follows from the classical Akizuki-Nakano vanishing theorem since $0+p<n$ and ${\cal O}_X(k-pd)$ is negative. If $k-pd \geq 0$, the holomorphic line bundle $-K_X + {\cal O}_X(k-pd)$ is ample (because $-K_X>0$), so writing ${\cal O}_X(k-pd) = K_X + (-K_X + {\cal O}_X(k-pd))$ we get $H^{0,\,p}(X,\,{\cal O}_X(k-pd))\simeq H^{n,\,p}(X,\,-K_X +{\cal O}_X(k-pd)) = 0$, where the last equality follows from the classical Kodaira vanishing theorem since $n+p >n$ and $-K_X + {\cal O}_X(k-pd)$ is positive.

We conclude that if the vanishings (\ref{eqn:multiple-vanishings}) hold, then $H^0(X,\,\Omega^p_X(k)) = 0$ (which is precisely the stated conclusion (\ref{eqn:hypersurfaces-projective-spaces_coh-vanishing})) thanks to the above sequence of injections and the vanishing of the last vector space therein.

\vspace{1ex}

$\bullet$ {\it Step $2$.} We will now prove the vanishings (\ref{eqn:multiple-vanishings}) under the assumption $d\geq 3$.

We will reduce the proof of these vanishings to proving similar vanishings for cohomology on $\Proj^{n+1}$ via the short exact sequence: \begin{eqnarray}\label{eqn:restriction_exact-seq}0\longrightarrow {\cal O}_{\Proj^{n+1}}(-d)\longrightarrow {\cal O}_{\Proj^{n+1}} & \longrightarrow & {\cal O}_X \longrightarrow 0 \\
 f & \longmapsto & f_{|X} \end{eqnarray} whose exactness follows from the equivalences: \begin{eqnarray*}f_{|X}=0 \iff \mbox{div}\, f\geq [X] \iff f\in{\cal O}_{\Proj^{n+1}}(-X) = {\cal O}_{\Proj^{n+1}}(-d).\end{eqnarray*}

Fix $i\in\{0,\dots, p-1\}$. Tensoring the short exact sequence (\ref{eqn:restriction_exact-seq}) by the holomorphic vector bundle $\Omega^{p-i}_{\Proj^{n+1}}(k-id)$, we get the short exact sequence: \begin{eqnarray*}\label{eqn:normal-sequence_X}0\longrightarrow \Omega^{p-i}_{\Proj^{n+1}}(k-(i+1)d)\longrightarrow \Omega^{p-i}_{\Proj^{n+1}}(k-id)\longrightarrow \Omega^{p-i}_{\Proj^{n+1}|X}(k-id)\longrightarrow 0.\end{eqnarray*} From the induced long exact sequence in cohomology, we extract \begin{eqnarray*}H^i(\Proj^{n+1},\,\Omega^{p-i}_{\Proj^{n+1}}(k-id))\longrightarrow H^i(X,\,\Omega^{p-i}_{\Proj^{n+1}|X}(k-id))\longrightarrow H^{i+1}(\Proj^{n+1},\,\Omega^{p-i}_{\Proj^{n+1}}(k-(i+1)d)).\end{eqnarray*}

This shows that, if the extreme terms $H^i(\Proj^{n+1},\,\Omega^{p-i}_{\Proj^{n+1}}(k-id))$ and $H^{i+1}(\Proj^{n+1},\,\Omega^{p-i}_{\Proj^{n+1}}(k-(i+1)d))$ vanish, then the middle term $H^i(X,\,\Omega^{p-i}_{\Proj^{n+1}|X}(k-id))$ vanishes as well, as needed in (\ref{eqn:multiple-vanishings}). Thus, it remains to prove the vanishing of the two extreme terms. 

We will use the following standard result (see e.g. Theorem 10.7 in [Dem97, chapter VII]) that completes the standard result cited above as (\ref{eqn:H_p0_O(k)}): \begin{eqnarray}\label{eqn:H_pq_O(k)} H^{p,\,q}(\Proj^{n+1},\,{\cal O}(k)) = 0\end{eqnarray} in each of the following cases: \begin{eqnarray*}(a) & & q\notin\{0,\,p,\,n+1\}; \\
    (b) & & q=0, \hspace{2ex} k\leq p \hspace{2ex} \mbox{and}\hspace{1ex} (k,\,p)\neq(0,\,0); \\
    (c) & & p=q\notin\{0,\,n+1\} \hspace{2ex} \mbox{and}\hspace{1ex} k\neq 0; \\
    (d) & &  q=n+1, \hspace{2ex} k\geq -n-1 + p \hspace{2ex} \mbox{and}\hspace{1ex} (k,\,p)\neq(0,\,n+1).     \end{eqnarray*}

  In our case, for the first extreme term above we get: \begin{eqnarray*}H^i(\Proj^{n+1},\,\Omega^{p-i}_{\Proj^{n+1}}(k-id))\simeq H^{p-i,\,i}(\Proj^{n+1},\,{\cal O}_{\Proj^{n+1}}(k-id))\end{eqnarray*} and this space vanishes automatically if $i\notin\{0,\,p-i,\,n+1\}$ by (a) in (\ref{eqn:H_pq_O(k)}). Meanwhile, if $i=0$, this space is $H^{p,\,0}(\Proj^{n+1},\,{\cal O}_{\Proj^{n+1}}(k)) = 0$ by (b) in (\ref{eqn:H_pq_O(k)}) because $k\leq p$. (Indeed, the last inequality is equivalent to $1 + (1-d)/2 \leq 0$, which is further equivalent to $d\geq 3$, which holds by hypothesis.) The case $i=p-i$ does not occur since $p$ is odd. Finally, if $i=n+1$, this space is $H^{p-n-1,\,n+1}(\Proj^{n+1},\,{\cal O}_{\Proj^{n+1}}(k-(n+1)d)) = 0$ because $p-n-1 = -p-2<0$.

  We conclude that, in all cases, the first extreme term above vanishes: $H^i(\Proj^{n+1},\,\Omega^{p-i}_{\Proj^{n+1}}(k-id)) = 0$ thanks to the hypothesis $d\geq 3$.

  As for the second extreme term above, we get: \begin{eqnarray*}H^{i+1}(\Proj^{n+1},\,\Omega^{p-i}_{\Proj^{n+1}}(k-(i+1)d))\simeq H^{p-i,\,i+1}(\Proj^{n+1},\,{\cal O}_{\Proj^{n+1}}(k-(i+1)d))\end{eqnarray*} and this space vanishes automatically if $i\notin\{-1,\,p-i-1,\,n\}$ by (a) in (\ref{eqn:H_pq_O(k)}). Meanwhile, if $i=-1$, this space is $H^{p+1,\,0}(\Proj^{n+1},\,{\cal O}_{\Proj^{n+1}}(k)) = 0$ by (b) in (\ref{eqn:H_pq_O(k)}) because $k\leq p+1$. (Indeed, the last inequality is equivalent to $(1-d)/2 \leq 0$, which is further equivalent to $d\geq 1$, which always holds.) If $i=p-i-1$, our space is $H^{i+1,\,i+1}(\Proj^{n+1},\,{\cal O}_{\Proj^{n+1}}(k-(i+1)d)) = 0$ by (c) in (\ref{eqn:H_pq_O(k)}) because $k-(i+1)d\neq 0$. (Indeed, the last inequality is equivalent to each of the following equivalent inequalities: $p+1 + (1-d)/2 \neq (p-i)d \iff p(d-2) + 2d\neq 3$ and the last one is always true because $p\geq 1$ and $d\geq 3$.) Finally, if $i=n$, our space is $H^{p-n,\,n+1}(\Proj^{n+1},\,{\cal O}_{\Proj^{n+1}}(k-(n+1)d)) = 0$ since $p-n<0$.

  We conclude that, in all cases, thanks to the hypothesis $d\geq 3$, the second extreme term above vanishes: $H^{i+1}(\Proj^{n+1},\,\Omega^{p-i}_{\Proj^{n+1}}(k-(i+1)d)) = 0$.

  The proof is complete.   \hfill $\Box$

\vspace{6ex}

\noindent {\bf References.} \\

\noindent [Ati71]\, M. Atiyah --- {\it Riemann Surfaces and Spin Structures} --- Ann. Sci. ENS {\bf 4}, no. 1 (1971), 47-62.

\vspace{1ex}

\noindent [B\"ar11]\, C. B\"ar --- {\it Spin Geometry} --- \url{https://www.math.uni-potsdam.de/fileadmin/user_upload/Prof-Geometrie/Dokumente/Lehre/Veranstaltungen/SS11/spingeo.pdf}

\vspace{1ex}

\noindent [Bea98]\, A. Beauville --- {\it Fano Contact Manifolds and Nilpotent Orbits} --- Comm. Math. Helv., {\bf 73} (4) (1998), 566--583.



\vspace{1ex}

\noindent [BDPP13]\, S. Boucksom, J.-P. Demailly, M. Paun, T. Peternell --- {\it The Pseudo-effective Cone of a Compact K\"ahler Manifold and Varieties of Negative Kodaira Dimension} --- J. Alg. Geom. {\bf 22} (2013) 201-248.

\vspace{1ex}

\noindent [Dem 84]\, J.-P. Demailly --- {\it Sur l'identit\'e de Bochner-Kodaira-Nakano en g\'eom\'etrie hermitienne} --- S\'eminaire d'analyse P. Lelong, P. Dolbeault, H. Skoda (editors) 1983/1984, Lecture Notes in Math., no. {\bf 1198}, Springer Verlag (1986), 88-97.

\vspace{1ex}

\noindent [Dem02]\, J.-P. Demailly --- {\it On the Frobenius Integrability of Certain Holomorphic $p$-Forms } --- In: Bauer, I., Catanese, F., Peternell, T., Kawamata, Y., Siu, YT. (eds) Complex Geometry. Springer, Berlin, Heidelberg. \url{https://doi.org/10.1007/978-3-642-56202-0_6}

\vspace{1ex}

\noindent [Dem97]\, J.-P. Demailly --- {\it Complex Analytic and Algebraic Geometry} --- \url{https://www-fourier.univ-grenoble-alpes.fr/~demailly/manuscripts/agbook.pdf}

\vspace{1ex}

\noindent [Die06]\, N. Q. Dieu --- {\it $q$-Plurisubharmonicity and $q$-Pseudoconvexity in $\C^n$} --- Publ. Mat., Barc. {\bf 50} (2006), no. 2, p. 349-369.

\vspace{1ex}

\noindent [Din22]\, S. Dinew --- {\it $m$-subharmonic and $m$-plurisubharmpnic functions: on two problems of Sadullaev} --- Ann. Fac. Sci. Toulouse Math. (6) {\bf 31} (2022), no. 3, 995-1009.

\vspace{1ex}

\noindent [DP25a]\, S. Dinew, D. Popovici --- {\it $m$-Positivity and Regularisation} --- arXiv:2510.25639v1 [math.DG].

\vspace{1ex}

\noindent [DP25b]\, S. Dinew, D. Popovici --- {\it $m$-Pseudo-effectivity and a Monge-Amp\`ere-Type Equation for Forms of Positive Degree} --- arXiv:2510.27362v1 [math.DG].

\vspace{1ex}

\noindent [HL13]\, F.R. Harvey, H.B. Lawson --- {\it p-Convexity, p-Plurisubharmonicity and the Levi Problem} --- Indiana Univ. Math. J. {\bf62} (2013), no. 1, 149-169.

\vspace{1ex}

\noindent [KPU25a]\, H. Kasuya, D. Popovici, L. Ugarte --- {\it Higher-Degree Holomorphic Contact Structures} ---  arXiv:2502.01447v2 [math.DG].

\vspace{1ex}

\noindent [KPU25b]\, H. Kasuya, D. Popovici, L. Ugarte --- {\it Properties of Holomorphic p-Contact Manifolds} --- arXiv:2511.10818v1 [math.DG].

\vspace{1ex}

\noindent [LM89]\, H. Blaine Lawson, M.-L. Michelsohn --- {\it Spin Geometry} --- Princeton University Press, Princeton, New Jersey, 1989.

\vspace{1ex}

\noindent [LeB95]\, C. Le Brun --- {\it Fano Manifolds, Contact Structures, and Quaternionic Geometry} --- Int. J. of Math. {\bf 6} (1995), 419–437.

\vspace{1ex}

\noindent [Lic63]\, A. Lichnerowicz --- {\it Spineurs harmoniques} --- C. R. Acad. Sci. Paris {\bf 257} (1963), 7–9.

\vspace{1ex}

\noindent [Ver10]\, M. Verbitsky --- {\it Plurisubharmonic Functions in Calibrated Geometry and q-Convexity} --- Math. Z. {\bf 264} (2010), no. 4, p. 939-957.

\vspace{1ex}

\noindent [TW10]\, V. Tosatti, B. Weinkove --- {\it The Complex Monge-Amp\`ere Equation on Compact Hermitian Manifolds} --- J. Amer. Math. Soc. {\bf 23} (2010), no. 4, 1187-1195.

\vspace{1ex}

\noindent [Yau78]\, S.T. Yau --- {\it On the Ricci Curvature of a Complex K\"ahler Manifold and the Complex Monge-Amp\`ere Equation I} --- Comm. Pure Appl. Math. {\bf 31} (1978) 339-411.

\vspace{6ex}

\noindent School of Maths and Physics     \hfill Institut de Math\'ematiques de Toulouse,

\noindent The University of Queensland     \hfill  Universit\'e de Toulouse,

\noindent St. Lucia, QLD 4067, Australia    \hfill  118 route de Narbonne, 31062 Toulouse, France

\noindent  Email: k.broder@uq.edu.au       \hfill     Email: popovici@math.univ-toulouse.fr

\end{document}